\documentclass[onefignum,onetabnum]{siamonline190516}
\usepackage{multirow}
\usepackage{lipsum}
\usepackage{epstopdf}
\usepackage{amsfonts,amsmath,amstext,amssymb,bm, amsopn, algorithmic,float,mathrsfs}
\usepackage{graphicx}
\usepackage{booktabs}
\usepackage{enumitem,soul}
\usepackage{multirow}
\usepackage[figuresright]{rotating}
\usepackage{siunitx}
\usepackage{array}
\usepackage{tabularx}

\newsiamremark{remark}{Remark}
\newsiamthm{example}{Example}

\makeatletter
\newcommand*{\addFileDependency}[1]{
  \typeout{(#1)}
  \@addtofilelist{#1}
  \IfFileExists{#1}{}{\typeout{No file #1.}}
}
\makeatother

\headers{Constrained GP for Parameter Inference}{Zhaohui. Li, Shihao. Yang, and C. F. Jeff. Wu}
\title{Parameter Inference based on Gaussian Processes Informed by Nonlinear Partial Differential Equations}
\author{Zhaohui Li\footnote{Joint first authors. These authors contributed equally.} $^{,}$\thanks{Academy of Mathematics and Systems Science, Chinese Academy of Sciences. 
School of Mathematical Sciences, University of Chinese Academy of Sciences, Beijing, China (\email{lizh@amss.ac.cn},).}
\and Shihao Yang\footnotemark[1] $^{,}$\thanks{H. Milton Stewart School of Industrial and Systems Engineering, Georgia Institute of Technology, Atlanta, GA, USA
  (\email{shihao.yang@isye.gatech.edu}, \email{jeff.wu@isye.gatech.edu}).}
\and C. F. Jeff Wu\footnotemark[3] $^{,}$\thanks{Corresponding Author. \funding{Wu’s research is supported by NSF DMS-1914632. Yang’s research is supported by NSF DMS-2318883.}}}

\ifpdf
\hypersetup{
  pdftitle={PIGPI},
  pdfauthor={Zhaohui Li}
}
\fi
\usepackage{soul}
\newcommand{\blue}{\textcolor{black}}

\begin{document}

\maketitle
\begin{abstract}
Partial differential equations (PDEs) are widely used for the description of physical and engineering phenomena.
Some key parameters involved in PDEs, which represent certain physical properties with important scientific interpretations, are difficult or even impossible to measure directly. 
Estimating these parameters from noisy and sparse experimental data of related physical quantities is an important task.
Many methods for PDE parameter inference involve a large number of evaluations for numerical solutions to PDE through algorithms such as the finite element method, which can be time-consuming, especially for nonlinear PDEs. 
In this paper, we propose a novel method for the inference of unknown parameters in PDEs, called the PDE-Informed Gaussian Process (PIGP) based parameter inference method.
Through modeling the PDE solution as a Gaussian process (GP), we derive the manifold constraints induced by the (linear) PDE structure such that, under the constraints, the GP satisfies the PDE.
\blue{For nonlinear PDEs, we propose an augmentation method that transforms the nonlinear PDE into an equivalent PDE system linear in all derivatives, which our PIGP-based method can handle. 
The proposed method can be applied to a broad spectrum of nonlinear PDEs. }
The PIGP-based method can be applied to multi-dimensional PDE systems and PDE systems with unobserved components. 
\blue{Like conventional Bayesian approaches,} the method can provide uncertainty quantification for both the unknown parameters and the PDE solution.
\blue{The PIGP-based method also completely bypasses the numerical solver for PDEs.} 
The proposed method is demonstrated through several application examples from different areas. 
\end{abstract}

\begin{keywords}
Gaussian Process, Partial Differential Equations, Model Calibration, Inverse Problems
\end{keywords}

\section{Introduction}\label{sec:introd}
Partial differential equations (PDEs) provide important theoretical understanding and computational implementation for engineering and physical systems. 
In many applications that involve PDE models, estimating/inferring key parameters, such as thermal conductivity in heat equations \cite{ozisik2021inverse}, hydraulic conductivity in groundwater flow equations \cite{li2005geostatistical}, from physical phenomenon data is a key task to ensure the PDE model accurately approximate the real-world phenomenon. 
PDE parameter inference is also known as inverse problems (e.g., \cite{kaipio2006statistical,mueller2012linear}) in many applications. 
One example is the estimation of the thermal conductivity of a material from the limited experimental data on temperature through heat PDE, which is essential for predicting its thermal behavior in various engineering applications \cite{taler2006solving}. 
Other examples include the inverse heat source problem \cite{yan2011new} and hydraulic transmissivity estimation problem \cite{li2005geostatistical}, etc. \par
\subsection{Literature Review}
Statistical methods \cite{muller2002fitting}, especially the Bayesian framework, which is also called Bayesian model calibration \cite{kennedy2001bayesian} in some statistical literature, have become more and more popular for PDE parameter inference.
Examples can be found in \cite{chiachio2021bayesian}. 
In the commonly used Bayesian framework, the evaluation of posterior density involves solving the PDE model with specified parameters. 
In practice, most PDE models are solved numerically using time-consuming algorithms such as finite element methods (FEM)\cite{johnson2012numerical} or finite difference methods (FDM)\cite{smith1985numerical}.
For direct inference using a Bayesian approach, repeated evaluations of the PDE solution are necessary, thereby exacerbating the computational burden of PDE parameter inference. 
To mitigate the computational burden caused by repeated evaluations of the numerical solution of PDE, there are generally two types of methods.\par
In the first type, a computer model emulation, also called a surrogate model, is established to approximate the PDE solution\cite{stuart2018posterior} or to estimate the likelihood and its transformations\cite{stuart2010inverse}.
For example, \cite{marzouk2007stochasticspectral} proposes a stochastic spectral method to solve the Bayesian inverse problem. 
\cite{lan2016emulation} employs the Gaussian process (GP) to approximate the PDE solution and its derivative while using Hessian to construct the Markov chain Monte Carlo algorithm. 
\cite{li2014adaptive} propose to use polynomial chaos to approximate the PDE solutions and an adaptive algorithm is proposed to sequentially construct the polynomial chaos surrogate model.
However, these surrogate-based methods for PDE parameter inference still suffer a serious computational burden (especially for nonlinear transient PDEs) since the construction of the surrogate model relies on a number of PDE numerical solutions\cite{santner2019design}. 
\blue{On the other hand, most surrogate models treat the numerical solver that produces the PDE solution as a black box. 
Thus the accurate construction of a surrogate model highly relies on the quality of experimental design \cite{sacks1989design}, which is a challenging task when parameter space is in high dimension.
Moreover, it is often challenging to ensure the extrapolation properties of the surrogate model \cite{queipo2005surrogate}.}

In the second type, numerical methods for solving PDEs are completely bypassed.
Mathematical models for function approximation, such as a linear basis function approximation or a GP regression\cite{rai2019gaussian}, are employed to approximate the PDE solution.

The simplest and most widely used approach that completely bypasses the numerical solver is called the two-stage method \cite{bar1999fitting,muller2004parameter,rai2019gaussian}.
In these methods, the PDE solution and its derivatives are reconstructed from observation data using basis function expansion, such as polynomial basis\cite{franke1998solving}, splines \cite{bar1999fitting}, etc, and coefficients for basis representation are estimated via regression methods.
Then the parameter $\bm{\theta}$ is estimated using the plug-in estimates of the PDE solution and its derivatives.
Although two-stage methods are simple to implement and efficient in computation, they have difficulties in estimating the PDE solution and its derivatives accurately, especially when the data is sparse and noisy\cite{muller2004parameter}. 
This will cause significant parameter estimation errors.

To mitigate this issue with the two-stage method when the observation data is sparse and noisy, a few \textit{integrated} approaches have been proposed.
Xun et.al., \cite{xun2013parameter} employs basis function expansion to approximate the PDE solution. 
Then a parameter cascading method and a Bayesian approach are proposed to simultaneously infer the coefficients of the basis function and parameters. 
The method is applied to linear PDEs with physical observations on a dense mesh grid. 
\cite{zhou2020inferring} replaces the basis function approximation in \cite{xun2013parameter} with the GP regression. 
A similar framework with the parameter cascading method in \cite{xun2013parameter} is proposed for linear PDE parameter inference. 
Moreover, they attempt to estimate solutions for nonlinear PDEs using Picard iterative linearization. 
However, the paper does not give a general methodology for solving parameter inference for nonlinear PDE. 
Meanwhile, this method does not quantify the uncertainty of the PDE solution that is proposed by the GP prediction. 
\blue{Employing a GP as a prior model for PDE solution has been extensively investigated\cite{graepel2003solving,cockayne2017probabilistic,raissi2017machine,chen2022apik,chkrebtii2016bayesian}.
For example, \cite{cockayne2017probabilistic} propose to use GPes as a prior model for PDE solution to solve forward and inverse problems.} 
The PDE solution is approximated by the GP conditional on the PDE constraint at collocation points. 
The posterior density of the unknown parameter is then obtained by replacing the PDE solution with the GP model.
The paper's limitation is that it attempts to solve inverse problems involving nonlinear PDEs, but it can only be applied to monotonic zeroth-order-nonlinear terms.
Constrained GP have garnered a lot of attention over the last decade \cite{swiler2020survey}.
\cite{raissi2017machine}, \cite{spitieris2023bayesian}, and \cite{oates2019bayesian} propose to use GP as a prior and apply the constrained GP model to solve the forward and inverse problems involving \textit{linear} PDEs. 
\cite{chen2021solving} assign a GP prior to the PDE solution. 
The inverse problem can be transformed to minimize the RKHS norm with constraints obtained by the PDE operator discretized on a collocation point set. 
The parameter is estimated by solving the optimization problem using the Gauss-Newton-type algorithm.

Another research topic that is relevant to PDE parameter inference is PDE identification \cite{liu2021automated,rudy2017data,chen2021gaussian}.
In this topic, the explicit form of PDE is unknown and must be inferred using sparse regression techniques.
PDE identification methods \cite{liu2021automated,rudy2017data,chen2021gaussian} are only applicable to identify the PDE that is linear with respect to the parameters.\par
Thanks to the vast development of deep learning in past decades, the Physics-Informed Neural Networks (PINN) \cite{raissi2019physics} have been proposed to solve forward problems, inverse problems \cite{depina2022application}, and PDE identification problems \cite{raissi2017physics} involve complex PDEs. 
However, using PINN for inverse problems poses challenges in computational efficiency when training a complex neural network.
Furthermore, how deep learning methods, including PINNs, handle \textit{uncertainty quantification} for parameters and PDE solutions remains unclear.. 

\subsection{Our Works and Contributions}
In this paper, we propose a novel method for PDE parameter inference, called \textbf{P}DE-\textbf{I}nformed \textbf{G}aussian \textbf{P}rocess (PIGP) based parameter inference method. 
The new method does not require any time-consuming numerical algorithms to solve the PDE and achieves simultaneous inference of unknown parameters and unobserved PDE solutions(under the inferred parameters) in one \textit{single} integrated step.
In particular, we assign a GP prior to the solution of PDE, inspired by the recent work in ODE space \cite{yang2021inference}. 
By restricting the GP on a manifold that satisfies the PDE, and their boundary and initial conditions, we are able to construct the posterior density for parameters and PDE solution.
For nonlinear PDE operators, we propose an augmentation method that constructs an equivalent PDE system with zeroth order nonlinearity, i.e., the nonlinear term of PDE only depends on the zeroth order derivatives of the PDE solution. (See Section \ref{subsec:PIGPI2} for details). 
As such, the time-consuming numerical PDE solver is completely bypassed and the computational burden is drastically reduced. \blue{In fact, our method exhibits scaling that is sub-quadratic, marginally above linear, with respect to the number of discretization points. This scaling efficiency depends on the smoothness of the PDE (See Section \ref{subsec:dimension_reduction} for details).}
Moreover, the proposed method provides uncertainty quantification for the parameter and the PDE solution simultaneously. 

\textbf{The main contributions of this paper} are summarised as follows. 
First, we propose a PIGP method for PDE parameter inference that completely bypasses the numerical solver. 
Second, in order to incorporate the nonlinear PDE operators, we propose a PDE augmentation method that constructs a zeroth-order-nonlinear PDE that is equivalent to the original PDE.
Through augmentation, we allow the PIGP method to handle a broad range of non-linear PDEs. 
Third, we propose a PIGP component that incorporates initial and boundary conditions, which are often available in practice. 
Hamiltonian Monte Carlo algorithm is utilized to draw random samples from the posterior density. 
Moreover, we propose to use the maximum a posteriori (MAP) estimation as an even faster summary of parameter inference, and a normal approximation for posterior density as an efficient approximated uncertainty quantification for the unknown parameters and PDE solution. 
Four examples are employed to demonstrate the performance of the proposed method.\par 
\textbf{Distinction of PIGP method from other GP-based inference methods:} Our work is not the first that employs a GP prior on the PDE solution for parameter inference. 
In general, there are two types of GP models employed in PDE parameter inference. 
The first type employs GP as a surrogate model to approximate the parameter-PDE solution relationship\cite{lan2016emulation}. 
These methods rely on the numerical solution of a PDE on a design as training data. 
Thus, these methods might suffer computational burden when the PDE is time-consuming to solve.
The second type \cite{cockayne2017probabilistic,zhou2020inferring} employs GP as a prior for the PDE solution. 
In these works, the PDE is generally assumed to be linear. 
\cite{cockayne2017probabilistic} considers the GP conditional on the linear constraints induced by the linear PDE. 
Then, the conditional GP is treated as a prediction of the PDE solution. 
The posterior for $\bm{\theta}$ is constructed by replacing the PDE solution with the conditional GP. 
\cite{zhou2020inferring} employs a similar framework with the parameter cascading method proposed in \cite{xun2013parameter} that replaces the basis function representation with the GP regression. 
The posterior for $\bm{\theta}$ is constructed by replacing the PDE solution with the GP prediction, conditional on the PDE constants. 

These existing methods, however, often rely on the linearity of the PDE to obtain the posterior distribution, while for nonlinear PDE, only heuristics are given for very limited nonlinear cases. 
As a comparison, our method is designed for solving PDE parameter inference where PDE can be nonlinear. 

The remainder of this paper is organised as follows. 
In section \ref{sec:method}, we introduce the PDE parameter inference problem and propose the PIGP method for PDE parameter inference. 
Section \ref{sec:algorithm} introduces the algorithm for implementation of the PIGP method. 
Section \ref{sec:simulation_res} illustrates PIGP method with several practical examples. 
Finally, section \ref{sec:conclusion} concludes the paper. 

\section{Methodology}\label{sec:method}

We begin by formally framing the \textit{PDE parameter inference} problems in Section \ref{subsec:basic_formulation}, and then introduce the PIGP-based inference methodology for PDE parameter inference in the remaining subsections. 
\subsection {PDE Parameter Inference} \label{subsec:basic_formulation}
In this paper, we consider parameterized and nonlinear partial differential equations (PDEs) that can be described by a derivative operator $\mathcal{A}(u,\bm{x}, \bm{\theta})$, which we call a PDE operator in this paper, and a function $f$, i.e.,
\begin{equation}\label{eq:general_pde}
	\mathcal{A}(u, \bm{x},\bm{\theta})=f(u, \bm{x},\bm{\theta}), \bm{x}\in \Omega 
\end{equation}
where $\bm{x}=(x_1,\dots,x_p)$, $u(\bm{x})$ is the solution of PDE defined on a domain $\Omega\subset \mathbb{R}^p$, 
$\mathcal{A}(u, \bm{x},\bm{\theta})$ is an $l$ dimensional differential operator on $u\in \mathcal{F}$ of order $a$, $\mathcal{F}$ is a Hilbert space. 
$f(u, \bm{x},\bm{\theta})$, which is called the zeroth-order term of PDE, is a $l$ dimensional vector-valued function of $u, \bm{x},\bm{\theta}$ but not any derivative of $u$. 
\par 
\begin{example}[LIDAR Equation]
Consider a long-range infrared light detection and ranging (LIDAR) equation \cite{xun2013parameter}, i.e.,
\begin{equation}\label{eq:lidar}
    \frac{\partial u(t,s)}{\partial t}-\theta_D \frac{\partial^2 u(t,s)}{\partial s^2}-\theta_S \frac{\partial u(t,s)}{\partial s} =\theta_A u(t,s).
\end{equation}
By definition, we have $\mathcal{A}(u, \bm{x},\bm{\theta}) = \frac{\partial u(t,s)}{\partial t}-\theta_D \frac{\partial^2 u(t,s)}{\partial s^2}-\theta_S \frac{\partial u(t,s)}{\partial s}$ and $f(\bm{x},u(\bm{x}),\bm{\theta}) = \theta_A u(t,s)$, where $\bm{x} = (t,s)$, $\bm{\theta} = (\theta_D,\theta_S,\theta_A)$ represents the unknown parameter. 
\end{example}
We assume that the solution of PDE is observed on a sparse set with measurement errors. 
Thus, we have observation data $(y_i,\bm{x}_i)$ where $y_i =u(\bm{x}_i)+\varepsilon_i, i=1,\dots,n$, where the random errors $\varepsilon_i, i=1,\dots,n$ are assumed independent and identically distributed (i.i.d.) normal random variables, i.e., $\varepsilon_i\sim N(0,\sigma_e^2), i=1,\dots,n$.
Let $\bm{\tau}=\{\bm{x}_1,\dots, \bm{x}_{n}\}$ denote the set of observation locations. 
The primary goal is to estimate $\bm{\theta}$ using observational data $y(\bm{\tau}) = \{y_i, i=1,\dots,n\}$ and quantify the uncertainty of our estimation.\par
In this paper, we propose a novel methodology, called PDE-Informed Gaussian Process (PIGP) based method, for estimating $\bm{\theta}$ and PDE solution $u(\bm{x})$ based on the data $y_{\bm{\tau}}$.
To infer the system parameters, we focus on the PDEs defined by (\ref{eq:general_pde}) where the zeroth-order term $f$ can be nonlinear, such as semi-linear PDEs (see examples for the linear zeroth-order term in Section \ref{subsec:toyeg1} and nonlinear zeroth-order term in Section \ref{subsec:egreaction_diffution}). 
Our method also accommodates nonlinear PDEs where the operator $\mathcal{A}$ depends on PDE solution $u$ \textit{non-linearly} (see examples in Section \ref{subsec:eg_burger_eq}). 
We begin with the introduction of the PIGP method for linear PDE operators in Section \ref{subsec:PIGPI1}. 
We will then discuss the more complicated non-linear PDEs in Section \ref{subsec:PIGPI2}.

\subsection{PIGP for PDE Parameter Inference with Linear PDE Operators} \label{subsec:PIGPI1}
We start with the linear PDE operators. 
To make a distinction, assume $\mathcal{L}_{\bm{x}}^{\bm{\theta}} u $ is linear with respect to $u$ in this section, which is a special case of the general PDE operator $\mathcal{A}(u, \bm{x},\bm{\theta})$.
The linear PDE operator is defined as follows. 
\begin{definition}(Linear PDE Operator)
A PDE operator $\mathcal{L}_{\bm{x}}^{\bm{\theta}}$ is called a \textbf{linear} PDE operator if it has the form
\begin{equation}\label{eq:PDE_operator_linear}
	\mathcal{L}_{\bm{x}}^{\bm{\theta}} u(\bm{x})=\sum_{\bm{\alpha}_i\in A} c_i(\bm{\theta},\bm{x})\nabla^{\bm{\alpha}_i}u(\bm{x}),
\end{equation}
where $\nabla^{\bm{\alpha}_i}u(\bm{x})= \frac{\partial^{|\bm{\alpha}_i|}u(\bm{x})}{\partial^{\alpha_{i1}}x_1\cdots \partial^{\alpha_{ip}}x_p}$ is a partial derivative operator, $\bm{\alpha}_i=(\alpha_{i1},\dots,\alpha_{ip})$, $\alpha_{ij}=0,1,2,\dots$, and $|\bm{\alpha}_i|=\sum_{j=1}^p \alpha_{ij}$.
$A=\{\bm{\alpha}_i,i=1,...,l\}$ is the index set of derivative orders. 
$c_i(\bm{\theta},\bm{x})$ is a function of $\bm{x}$ and $\bm{\theta}$.
Moreover, $\alpha = \max_{\bm{\alpha}_i\in A}\{\|\bm{\alpha_i}\|\}$ is called the order of $\mathcal{L}_{\bm{x}}^{\bm{\theta}}$.
\end{definition}

We now introduce the fundamental concept behind the PIGP method for PDE parameter inference in cases where the PDE operator is linear. 
In general, PDE solution $u$ is assumed to be located in some Sobolev space $\mathcal{F}$ in which the proper degree of smoothness properties are satisfied such that the PDE operator is well defined. 
Thus following the Bayesian paradigm, we assume the PDE solution $u(\bm{x})$ is a realization of a GP $U(\bm{x})$.
In particular, assume that $U(\bm{x})\sim \text{GP}(\mu(\bm{x}),\mathcal{K}(\bm{x},\bm{x}'))$, where $\mu(\bm{x}),\mathcal{K}(\bm{x},\bm{x}')$ are the mean and covariance function of the GP. 
Let $\mathcal{H}$ denote the function space of all sample paths for the GP $U(\bm{x})$. 

To incorporate PDE information into the GP prior, we first define a manifold introduced by the PDE structure $\mathcal{U}=\{U| \mathcal{L}_{\bm{x}}^{\bm{\theta}} U(\bm{x}) = f(\bm{x},U(\bm{x}),\bm{\theta}), \theta\in\Omega_{\theta}\}$, where $\Omega_{\theta}$ is the parameter space.
In other words, $U\in \mathcal{U}$ lies on the manifold of the PDE solutions with some $\theta$ in parameter space. 
Then we introduce a new variable $W$ quantifying the difference between the GP sample path $U(\bm{x})$ and the PDE structure with a given value of the parameter $\bm{\theta}$, i.e., 
\begin{equation}
	W=\sup_{\bm{x}\in\Omega} \|\mathcal{L}_{\bm{x}}^{\bm{\theta}} U(\bm{x})-f(\bm{x},U(\bm{x}),\bm{\theta})\|_{\infty}. 
\end{equation}
Obviously, $W\equiv 0$ if and only if the PDE is satisfied by $U$ with given parameters $\bm{\theta}$, i.e., $U\in \mathcal{U}$.
The connection between $\mathcal{U}$ and $W$ justifies that conditioning on $W= 0$ would produce an ideal inference. 
However, $W$ is not computable in practice. 
Therefore, we approximate $W$ by a finite discretization on the set $\bm{I}=\{\bm{x}_1,\dots,\bm{x}_{n_I}\}\subset \Omega$ such that $\bm{\tau} \subset \bm{I} \subset\Omega$ and similarly define $W_{\bm{I}}$ as
\begin{equation}
	W_{\bm{I}}=\max_{\bm{x}\in\bm{I}} \|\mathcal{L}_{\bm{x}}^{\bm{\theta}} U(\bm{x})-f(\bm{x},U(\bm{x}),\bm{\theta})\|_{\infty}. 
\end{equation}
Intuitively, under mild continuity and smoothness conditions for $U$ and $f$, $W_{\bm{I}} \to W$ as $\bm{I}$ get dense in $\Omega$.
Thus, the constraint $W_{\bm{I}} = 0$ provides an approximation to $W=0$.

Based on the previous discussion, we acquiescently assume that $\mathcal{L}_{\bm{x}}^{\bm{\theta}} U(\bm{x})$ exists.
However, this is not a trivial assumption for an arbitrary GP model. 
To avoid ill-posedness of $\mathcal{L}_{\bm{x}}^{\bm{\theta}} U(\bm{x})$, we assume that the mean and kernel function of GP $U$ is chosen to be \textit{smooth enough} (differentiable up to at least $2a$ orders) such that the derivatives of $U$ and their distributions are all well defined. 
In particular, we herein assume the mean function is a constant for simplicity, i.e., $\mu(\bm{x} ) = \mu$.
Moreover, we assume that $\mathcal{K}$ is chosen that $\mathcal{L}_{\bm{x}}^{\bm{\theta}} \mathcal{K}(\bm{x},\bm{x}')$, $\mathcal{L}_{\bm{x}'}^{\bm{\theta}} \mathcal{K}(\bm{x},\bm{x}')$, and $\mathcal{L}_{\bm{x}}^{\bm{\theta}} \mathcal{L}_{\bm{x}'}^{\bm{\theta}} \mathcal{K}(\bm{x},\bm{x}')$ exist. 
Note that $\mathcal{L}_{\bm{x}}^{\bm{\theta}} U(\bm{x})$ is also GP conditional on $U(\bm{x})$ and $\bm{\theta}$ due to the \textit{linearity} of $\mathcal{L}_{\bm{x}}^{\bm{\theta}}$ (see \cite{adler2010geometry}, Theorem 2.2.2).

For notational simplicity, we use the notation $u(\bm{I}) = (u(\bm{x}_1),u(\bm{x}_2),\dots,u(\bm{x}_{n_{\bm{I}}}))$ and similarly for $U(\bm{I})$.
Let $f(\bm{I},u(\bm{I}),\bm{\theta}) = (f(\bm{x}_1,u(\bm{x}_1),\bm{\theta}),f(\bm{x}_2,u(\bm{x}_2),\bm{\theta}),\dots,f(\bm{x}_{n_{\bm{I}}},u(\bm{x}_{n_{\bm{I}}}),\bm{\theta}))$.
The posterior density is given by 
\begin{align}\label{eq:general_posterior}
	&p\left(\sigma_e^2,\bm{\theta}, u(\bm{I})|W_{\bm{I}}=0, Y(\bm{\tau})=y(\bm{\tau})\right)\notag\\
	\propto &\pi(\sigma_e^2)\times\pi_{\bm{\Theta}}\left(\bm{\theta}\right)\times P\left(U(\bm{I})=u(\bm{I})|\bm{\Theta}=\bm{\theta}\right)\notag\\
	&\times P\left(Y(\bm{\tau})=y(\bm{\tau})|\sigma_e^2,U(\bm{I})=u(\bm{I}), \bm{\Theta}=\bm{\theta}\right)\notag\\
	&\times P\left(W_{\bm{I}}=0|Y(\bm{\tau})=y(\bm{\tau}), U(\bm{I})=u(\bm{I}), \bm{\Theta}=\bm{\theta}\right)\notag\\
	= &\pi(\sigma_e^2)\times\pi_{\bm{\Theta}}\left(\bm{\theta}\right)\times \exp \Big\{-\frac{1}{2}\big[n_{\bm{I}}\log(2\pi) + \log(|C|)+\left\| u(\bm{I})-\mu(\bm{I})\right\|_{C^{-1}}\notag\\ 
	&+n\log(2\pi) + n \log(\sigma_e^2)+\left\| u(\bm{\tau})-y(\bm{\tau})\right\|_{\sigma_e^{-2}}\notag\\ 
	& + n_{\bm{I}}\log(2\pi) +\log(|K|)+\left\|f(\bm{I},u(\bm{I}),\bm{\theta})-\mathcal{L}_{\bm{x}}^{\bm{\theta}}\mu(\bm{I})-m\{u(\bm{I})-\mu(\bm{I})\} \right\|_{K^{-1}}\big]  \Big\},
\end{align}
where
\begin{equation*}
	\left\{
	\begin{aligned}
		C&=\mathcal{K}(\bm{I},\bm{I})\\
		m&=\mathcal{LK}(\bm{I},\bm{I})\mathcal{K}(\bm{I},\bm{I})^{-1}\\
		K&=\mathcal{LKL}(\bm{I},\bm{I})-\mathcal{LK}(\bm{I},\bm{I})\mathcal{K}(\bm{I},\bm{I})^{-1}\mathcal{KL}(\bm{I},\bm{I})
	\end{aligned}
	\right.,
\end{equation*}
where $\mathcal{K}(\bm{I},\bm{I})$ denotes an $N\times N$ matrix with $(i,j)$ element $\mathcal{K}(\bm{x}_i,\bm{x}_j)$ and similarly for $\mathcal{LK}, \mathcal{KL}, $ and $\mathcal{LKL}$;
$\mathcal{LK}(\bm{x},\bm{x}')=\mathcal{L}_{\bm{x}}^{\bm{\theta}}(\mathcal{K}(\bm{x},\bm{x}'))$,
$\mathcal{LKL}(\bm{x},\bm{x}')=\mathcal{L}_{\bm{x}}^{\bm{\theta}}(\mathcal{L}_{\bm{x}'}^{\bm{\theta}}(\mathcal{K}(\bm{x},\bm{x}')))$, and  $\mathcal{KL}(\bm{x},\bm{x}')=\mathcal{L}_{\bm{x}'}^{\bm{\theta}}(\mathcal{K}(\bm{x},\bm{x}'))$.
Appendix \ref{spsubsec:matern} has details for calculating these quantities. 
For the determination of the priors for $\bm{\theta}$ and $\sigma_e^2$, informative priors are encouraged if they are available from historical data or other resources.
Otherwise, uniform prior for $\bm{\theta}$ and Jeffrey's prior for $\sigma_e^2$ are chosen as a default non-informative prior.

The derivative of \eqref{eq:general_posterior} is simple. 
First, $P\left(U(\bm{I})=u(\bm{I})|\bm{\Theta}=\bm{\theta}\right)=P\left(U(\bm{I})=u(\bm{I})\right)$ is Gaussian due to the GP prior for $U(\bm{x})$;
Second, $P\left(Y(\bm{\tau})=y(\bm{\tau})|\sigma_e^2, U(\bm{I}) =u(\bm{I}), \bm{\Theta}=\bm{\theta}\right)$ is Gaussian due to the distribution of random measurement errors;
Third, 
\begin{align*}
	&P\left(W_{\bm{I}}=0|Y(\bm{\tau})=y(\bm{\tau}), U(\bm{I})=u(\bm{I}), \bm{\Theta}=\bm{\theta}\right)\\
	=& P\left(\mathcal{L}_{\bm{x}}^{\bm{\theta}}U(\bm{I})= f(\bm{I},u(\bm{I}),\bm{\theta})| Y(\bm{\tau})=y(\bm{\tau}), U(\bm{I})=u(\bm{I}), \bm{\Theta}=\bm{\theta} \right)\\
	=& P\left(\mathcal{L}_{\bm{x}}^{\bm{\theta}}U(\bm{I})= f(\bm{I},u(\bm{I}),\bm{\theta})|U(\bm{I})=u(\bm{I}),  \bm{\Theta}=\bm{\theta} \right),
\end{align*}
is Gaussian due to the Gaussianity of $\mathcal{L}_{\bm{x}}^{\bm{\theta}}U(\bm{I})$ conditional on $U(\bm{I}) = u(\bm{I})$, the density is given in the last line of \eqref{eq:general_posterior}.\par
\begin{remark}
It is worth noting that our method can be directly applied to the problem where the observation is given by the form $y_i =h(u(\bm{x}_i))+\varepsilon_i, i=1,\dots,n$, where $h: \mathbb{R}\to \mathbb{R}$ is a known one-to-one map in $\mathcal{R}$, the random errors $\varepsilon_i, i=1,\dots,n$ are i.i.d. normal random variables, i.e., $\varepsilon_i\sim N(0,\sigma_e^2), i=1,\dots,n$.
To this end, one only needs to modify \eqref{eq:general_posterior} by replacing $\left\| u(\bm{\tau})-y(\bm{\tau})\right\|_{\sigma_e^{-2}}$ with $\left\| h(u(\bm{\tau}))-y(\bm{\tau})\right\|_{\sigma_e^{-2}}$. 
\end{remark}

In this paper, we employ the product Mat\'{e}rn kernel as the correlation function for the GP prior. 
The main reason for this choice is its ability to tune the variability of the GP in each direction.
In particular, the product Mat\'{e}rn kernel is given by $\mathcal{K}(\bm{x},\bm{x}')=\prod_{i=1}^p\phi_{1i} \frac{2^{1-\nu}}{\Gamma(\nu)}(\sqrt{2\nu}\frac{d_i}{\phi_{2i}})^{\nu}B_{\nu}(\sqrt{2\nu}\frac{d_i}{\phi_{2i}})$, where $d_i=|x_i-x_i'|, i=1,\dots, p$, $\Gamma$ is the Gamma function, $B_{\nu}$ is the modified Bessel function of the second kind, and the degree of freedom $\nu$ is set to be $2a+\delta$ to ensure that the $2a$-th order derivatives of the kernel with respect to any coordinate $x_1$ exists, where $\delta$ is a small positive number.
In this paper we take $\delta = 0.1$ and a universal degree of freedom $a = \max_{\bm{\alpha}_i\in A}\{\|\bm{\alpha_i}\|_{\infty}\}$ for all dimensions of $\bm{x}$.
Appendix \ref{spsubsec:kernel_choice} has detailed discussion on the choice of kernel function.\par
\begin{remark}
	In this paper, we estimate the GP hyper-parameters using only the observation data $\bm{y}_{\bm{\tau}}$ through maximum marginal likelihood and then we fix the hyper-parameters at their optimized values in posterior inference. 
	This is also called modularization and is justified in \cite{bayarri2009modularization}.
	Appendix \ref{spsubsec:train_GP} has details on the training of GP. 
\end{remark} 
\begin{remark}\label{rmk:missing_comp}
For notational simplicity, we assume PDE is a one-dimensional function in this section.
However, it is natural to generalize the PDE model (\ref{eq:general_pde}) to a multi-dimensional output PDE system with a multi-dimensional PDE solution. 
Examples can be found in Section \ref{subsec:egreaction_diffution}. 
Furthermore, our method can naturally handle the situation where one or more components in multi-dimensional output PDE systems are (partially) missing. 
We demonstrate the performance of PIGP method on PDE parameter inference based on censored data in Section \ref{subsec:egreaction_diffution}.
The capacity to handle missing components is crucial for incorporating general nonlinear PDE cases.
More discussion can be found in Section \ref{subsec:PIGPI2}. 
\end{remark}

\subsection{Incorporate Nonlinear PDE Operators}\label{subsec:PIGPI2}
Two limitations arise when applying PIGP to parameter inference with general PDEs:
\begin{itemize}
    \item \textbf{Computationally Expensive for Parameter Dependent Operator:} When $\mathcal{L}_{\bm{x}}^{\bm{\theta}}$ is linear with respect to $u$ but depends on parameter $\bm{\theta}$, the updating of $\bm{\theta}$ requires the updating of $\mathcal{LKL}$ and $K^{-1}$. This can considerably increase the computational complexity of posterior inference, especially for MCMC algorithms.
    \item \textbf{Non-Flexible for Non-Linear Operator:} For a non-linear operator $\mathcal{A}(u,\bm{x}, \bm{\theta})$ (refer to the nonlinear heat equation in Table \ref{tb:aug_example} as an example), $\mathcal{A}(U,\bm{x}, \bm{\theta})$ may not be Gaussian. 
    This prevents us from directly writing down the posterior density as given in \eqref{eq:general_posterior}.
\end{itemize}
Motivated by these limitations, we propose an augmentation step to overcome them.
We consider the nonlinear PDE defined in \eqref{eq:general_pde} of the form
\begin{equation}\label{eq:pde_operator_nonlinear}
  \mathcal{A}(u,\bm{x}, \bm{\theta})=f(u,\bm{x}, \bm{\theta}) \overset{\text{rewrite}}{\Longleftrightarrow}  \nabla^{\bm{\alpha}_1}u = \mathcal{A}_1(\bm{\theta},u,\nabla^{\bm{\alpha}_{2}}u,\dots,\nabla^{\bm{\alpha}_{l}}u) + f(u,\bm{x}, \bm{\theta})
\end{equation}
 where $\mathcal{A}_1$ is a nonlinear function of $(\bm{\theta},u,\nabla^{\bm{\alpha}_{2}}u,\dots,\nabla^{\bm{\alpha}_{l}}u)$, which is the remaining part of $\mathcal{A}$.
$\mathcal{A}_1$ may contain parameter-dependent components and nonlinear components. 
In this section, we construct an augmented PDE that has several advantages.
First, the augmented PDE is linear in all non-zeroth order derivatives.
Second, the PDE operator is independent of $\bm{\theta}$.
Third, the corresponding augmented PDE is equivalent to the original PDE.
To illustrate this, we introduce a nonlinear PDE example. 
\begin{example}[Burger's Equation]
Consider the Burger's equation,
\begin{equation}\label{eq:burgersequation}
    \frac{\partial u}{\partial t} - \theta_1 u \frac{\partial u}{\partial s} + \theta_2 \frac{\partial^2 u}{\partial s^2} = 0.
\end{equation}
Burger's equation is a nonlinear equation (due to the existence of $\theta_1 u \frac{\partial u}{\partial s}$) with parameters $\theta_1$ and $\theta_2$ involved in the PDE operator where $\mathcal{A}(u,\bm{x}, \bm{\theta}) u(\bm{x}) = \frac{\partial u}{\partial t} - \theta_1 u \frac{\partial u}{\partial s} + \theta_2 \frac{\partial^2 u}{\partial s^2}$ and $f = 0$.
Define an augmented PDE system $(u_1, u_2, u_3)$ with PDE constraints:
\begin{equation}\label{eq:PDE_augemented_burger}
\begin{cases}
    \frac{\partial u_1}{\partial s}(\bm{x}) &= u_2(\bm{x}),\\
    \frac{\partial u_2}{\partial s}(\bm{x})& = u_3(\bm{x}),\\
    \frac{\partial u_1}{\partial t} (\bm{x})&= \theta_1 u_1(\bm{x})u_2(\bm{x})-\theta_2 u_3 (\bm{x}).
\end{cases}
\end{equation}
Clearly, $u_1=u$ for all $\bm{x}\in \Omega$ if they are the solution of PDEs (\ref{eq:burgersequation}) and (\ref{eq:PDE_augemented_burger}) given $\bm{\theta}$, respectively.
We call (\ref{eq:PDE_augemented_burger}) an augmented PDE of (\ref{eq:burgersequation}). 
By assuming the augmented PDE given in (\ref{eq:PDE_augemented_burger}) and observation data given by $y_i = u_1(\bm{x}_i)+\varepsilon_i$, the PIGP method in Section \ref{subsec:PIGPI1} is immediately applicable. 
Note that for the augmentation PDE, we assume the observations for $u_2$ and $u_3$ are missing, which can be handled by PIGP as discussed in Remark \ref{rmk:missing_comp}. 
\end{example}
Now we provide a formal definition for the augmented PDE.
\begin{definition}[Augmented PDE]
Suppose a nonlinear PDE has the form $\mathcal{A}(u,\bm{x}, \bm{\theta}) = f(u,\bm{x}, \bm{\theta})$ where $\mathcal{A}(u,\bm{x}, \bm{\theta})$ is given in (\ref{eq:pde_operator_nonlinear}).
We call a PDE system $\bm{\mathcal{L}_x^{\mathcal{A}}}\bm{u}=\bm{f}^{\mathcal{A}}$ is an augmented PDE of $\mathcal{A}(u,\bm{x}, \bm{\theta}) = f(u,\bm{x}, \bm{\theta})$, if \\
(a). Linearity: $\bm{\mathcal{L}_x^{\mathcal{A}}}\bm{u}$ is a linear PDE operator, where $\bm{u} = (u_1,\dots,u_l)^T$;\\
(b). Invariance: $\bm{\mathcal{L}_x^{\mathcal{A}}}\bm{u}$ is independent to parameter $\bm{\theta}$;\\
(c).Equivalence: If $\bm{u}$ is a solution of $\bm{\mathcal{L}_x^{\mathcal{A}}}\bm{u}=\bm{f}^{\mathcal{A}}$, then $u=u_1$ is a solution of $\mathcal{A}(u,\bm{x}, \bm{\theta}) = f(u,\bm{x}, \bm{\theta})$. 
And vice versa, i.e., if $u$ is a solution of $\mathcal{A}(u,\bm{x}, \bm{\theta}) = f(u,\bm{x}, \bm{\theta})$, then $\bm{\mathcal{L}_x^{\mathcal{A}}}\bm{u}=\bm{f}^{\mathcal{A}}$ has a solution such that $u_1 = u$. 
\end{definition}
A method for constructing the augmented PDE is provided in this section.
Given the definition, the most straightforward way to construct the augmented PDE from the original PDE Eq.\eqref{eq:pde_operator_nonlinear} as follows.
We first rewrite \eqref{eq:pde_operator_nonlinear} as
\begin{equation*}
    \nabla^{\bm{\alpha}_1}u = \mathcal{A}_1(\bm{x},\bm{\theta},\underbrace{u}_{u_1},\underbrace{\nabla^{\bm{\alpha}_{2}}u}_{u_2},\dots,\underbrace{\nabla^{\bm{\alpha}_{l}}u}_{u_l}) + f(\bm{x},u,\bm{\theta}),
\end{equation*}
such that
\begin{equation}\label{eq:PDE_augemented}
    \begin{aligned}
    \nabla^{\bm{\alpha}_{2}} u_1(\bm{x})&= u_{2}(\bm{x}),\\
    ...\\
    \nabla^{\bm{\alpha}_{l}}u_1(\bm{x}) &= u_{l}(\bm{x}),\\
    \nabla^{\bm{\alpha}_{1}} u_1(\bm{x}) &= f(\bm{x},u_1(\bm{x}),\bm{\theta}) +\mathcal{A}_1(\bm{x},\bm{\theta},u_1,u_2,\dots,u_l).
\end{aligned}
\end{equation}
Due to the equivalence between augmented PDE (\ref{eq:PDE_augemented}) and (\ref{eq:pde_operator_nonlinear}), we propose to use the method in Section \ref{subsec:PIGPI1} and the PDE operator in (\ref{eq:PDE_augemented}) to generate the parameter inference framework.
In particular, we define the augmented operator $\bm{\mathcal{L}_x^{\mathcal{A}}}$ as a $l\times l$ matrix. 
For (\ref{eq:PDE_augemented}), the corresponding operator is given by
\begin{equation}
\bm{\mathcal{L}_x^{\mathcal{A}}} =\left[
\begin{aligned}
    \nabla^{\bm{\alpha}_2},0,\cdots,0\\
    \nabla^{\bm{\alpha}_3},0,\cdots,0\\
    \cdots\\
    \nabla^{\bm{\alpha}_1},0,\cdots,0
\end{aligned}
\right].
\end{equation}
We give some examples in Table \ref{tb:aug_example} to demonstrate the construction of augmented PDE from the original PDE. 
It can be seen that an augmented PDE is easily constructed for a large number of PDE models. 
\begin{table}[H]
\caption{Some examples for constructing PDE augmentation}\label{tb:aug_example}    
\begin{center}
\begin{tabular}{|p{4.2cm}|p{4.8cm}|p{5cm}|}
\hline
Name & Original Form & Augmented Form \\ \hline
Fisher's Equation  & 
$\frac{\partial u}{\partial t}-D \frac{\partial^{2} u}{\partial s^{2}}=r u(1-u)$ & 
\begin{tabular}{l}
    $\frac{\partial^{2} u_1}{\partial s^{2}} = u_2$ \\
    $\frac{\partial u_1}{\partial t} = D u_2+r u_1(1-u_1)$
\end{tabular}\\
\hline
Telegraph Equation & 
$\frac{\partial^{2} u}{\partial t^{2}}+k \frac{\partial u}{\partial t}=a^{2} \frac{\partial^{2} u}{\partial s^{2}}+b u$ & 
\begin{tabular}{l}
    $\frac{\partial u_1}{\partial t} = u_2$ \\
    $\frac{\partial^{2} u_1}{\partial s^{2}} = u_3$ \\
    $\frac{\partial u_2}{\partial t} = a^{2} u_3+b u_1 -k u_2$
\end{tabular} \\
\hline
Nonlinear Heat Equation & 
$\frac{\partial u}{\partial t}=a \frac{\partial}{\partial s}\left(e^{\lambda u} \frac{\partial u}{\partial s}\right)+b+c_{1} e^{\beta u}+c_{2} e^{\gamma u}$ & 
\begin{tabular}{l}
    $\frac{\partial u_1}{\partial s} = u_2$ \\
    $\frac{\partial u_2}{\partial s} = u_3$ \\
    $\frac{\partial u_1}{\partial t} = a\lambda e^{\lambda u_1} u_2^2 + a e^{\lambda u_1} u_3$\\ $+b+c_{1} e^{\beta u_1}+c_{2} e^{\gamma u_1}$
\end{tabular}  \\ 
\hline
Generalized Korteweg–de Vries Equation & 
$\frac{\partial u}{\partial t}+\frac{\partial^{3} u}{\partial s^{3}}+g(u) \frac{\partial u}{\partial s}=0$ &
\begin{tabular}{l}
    $\frac{\partial u_1}{\partial s} = u_2$ \\
    $\frac{\partial u_1}{\partial t}+\frac{\partial^{3} u_2}{\partial s^2} = g(u_1) u_2$
\end{tabular}  \\ 
\hline
Reaction-Diffusion System & 
\begin{tabular}{l}
$\frac{\partial u}{\partial t}=a \frac{\partial^2 u}{\partial s^2}+F(u, v)$,\\
$\frac{\partial v}{\partial t}=a \frac{\partial^2 v}{\partial s^2}+G(u, v)$ 
\end{tabular}
& 
\begin{tabular}{l}
    $\frac{\partial^2 u_1}{\partial s^2} = u_2$ \\
    $\frac{\partial^2 v_1}{\partial s^2} = v_2$ \\
    $\frac{\partial u_1}{\partial t} =a u_2+F(u_1, v_1)$\\
    $\frac{\partial v_1}{\partial t} =a v_2+G(u_1, v_1)$
\end{tabular}  \\ 
\hline
\end{tabular}    
\end{center}
\end{table}
PIGP method for augmented PDE is now quite straightforward.
We assign $U_1,U_2,U_3,\dots,U_l$ independent GP prior with mean and variance functions $\mu_i,\mathcal{K}_i, i=1,\dots,l$ respectively. 
The choice of $\mu_i,\mathcal{K}_i, i=1,\dots,l$ are discussed in Appendix \ref{spsubsec:kernel_choice}.
Now we can easily write down the posterior density of $u(\bm{I}),u_2(\bm{I}),\dots,u_l(\bm{I}),\bm{\theta}$. 
For notational simplicity, we write $\bm{u}(\bm{I})=(u_1(\bm{I}),u_2(\bm{I}),\dots,u_l(\bm{I}))$. 
The posterior can be written as the form (\ref{eq:general_posterior}) by replacing $\mathcal{L}$ and $f$ with $\bm{\mathcal{L}_x^{\mathcal{A}}}$ and $\bm{f}^{\mathcal{A}}$ respectively, where
\begin{equation*}
	\left\{
	\begin{aligned}
		&C = \text{diag}(\mathcal{K}_1(\bm{I},\bm{I}),\dots,\mathcal{K}_l(\bm{I},\bm{I}))\\
		&\bm{m}=\bm{\mathcal{L}_x^{\mathcal{A}}}\mathcal{K}(\bm{I},\bm{I})\mathcal{K}(\bm{I},\bm{I})^{-1} \\
		&K=\bm{\mathcal{L}_x^{\mathcal{A}}}\mathcal{K}\bm{\mathcal{L}_{x'}^{\mathcal{A}}}(\bm{I},\bm{I})-\bm{\mathcal{L}_x^{\mathcal{A}}}\mathcal{K}(\bm{I},\bm{I})\mathcal{K}(\bm{I},\bm{I})^{-1}\mathcal{K}\bm{\mathcal{L}_x^{\mathcal{A}}}(\bm{I},\bm{I})\\
		&\bm{\mathcal{L}_x^{\mathcal{A}}}\mathcal{K}(\bm{I},\bm{I}) = \bm{\mathcal{L}_x^{\mathcal{A}}}\mathcal{K}\otimes (\bm{I})\\
		& \bm{\mathcal{L}_x^{\mathcal{A}}}\mathcal{K}\bm{\mathcal{L}_{x'}^{\mathcal{A}}}(\bm{I},\bm{I}) =\bm{\mathcal{L}_x^{\mathcal{A}}}\mathcal{K}{\bm{\mathcal{L}_x^{\mathcal{A}}}}^T\otimes (\bm{I})
	\end{aligned}
	\right.,
\end{equation*}
where $\otimes$ is the tensor product, which defines an operation over a matrix-type operator and a point set, i.e., $\mathcal{P}\otimes(\bm{I}) =(\mathcal{P}_{ij} (\bm{I}))$ is a $ln_{\bm{I}}\times ln_{\bm{I}}$ matrix whose $i,j$ block is $\mathcal{P}_{ij} (\bm{I})$, $\mathcal{P}$ is an operator whose $i,j$ element is $\mathcal{P}_{ij}$.
Since we assume the prior of the components of $\bm{U}$ are independent, $\mathcal{K}$ is a block diagonal matrix-valued function whose diagonal elements are $\mathcal{K}_1,\mathcal{K}_2,\dots, \mathcal{K}_l$.
We discuss the training of GP regression in the supplementary materials, see Appendix \ref{spsubsec:train_GP}.\par
It is important to note that the augmentation is not unique. 
For a specific PDE, there might be more than one augmented PDE that exists and satisfies the definition. 
For example, another augmented PDE for (\ref{eq:burgersequation}) is 
\begin{equation}\label{eq:burger_apde2}
\begin{aligned}
\frac{\partial u_1(\bm{x})}{\partial s} &= u_2(\bm{x}),\\
\frac{\partial^2 u_1(\bm{x})}{\partial s^2}&= u_3(\bm{x}),\\
\frac{\partial u_1(\bm{x})}{\partial t} &= \theta_1 u_1(\bm{x})u_2(\bm{x})-\theta_2 u_3(\bm{x}).
\end{aligned}
\end{equation}
Our empirical experiments suggest that both methods work well, and our method doesn't rely heavily on the choice of augmented PDE. 
However, we can benefit from carefully choosing the augmented PDE from its multiple candidates. 
To clarify, we suggest a rule of thumb for selecting an augmented PDE. 
\begin{remark}
\textbf{Lowest order of derivative (LOD)} principal: 
As an example, we recommend using the augmented PDE (\ref{eq:PDE_augemented_burger}) instead of (\ref{eq:burger_apde2}) for two reasons.
First, (\ref{eq:PDE_augemented_burger}) is a first-order PDE while the latter is a second-order PDE.
We intend to choose an augmented PDE with \textit{lowest order of derivative}, which is helpful for improving the numerical stability of GP. 
Second, since we assign independent GP priors for $u_i$, $u_j$ where $i\neq j$, the covariance matrix $K$ is a block-wise matrix consisting of $l\times l$ components. 
The covariance matrix $K$ associated with (\ref{eq:PDE_augemented_burger}) has 4 zero blocks, which is helpful to improve computational efficiency.
In summary, when there is multiple augmented PDE available for a nonlinear PDE model, we recommend choosing one that has the lowest order of derivatives.  
\end{remark}
\begin{example}[Counterexample: Eikonal Equation]
At the end of this subsection, we give an example to which our framework is not applicable.
The PDE is a version of the well-known Eikonal equation \cite{gremaud2006computational},
\begin{equation}\label{counter_eg}
     (\frac{\partial u}{\partial x_1})^2 + (\frac{\partial u}{\partial x_2})^2 = f(u,\bm{x},\bm{\theta}),
\end{equation}
where $f(u,\bm{x},\bm{\theta})$ is a function with positive values.
By simple algebra, we get two PDEs, i.e., 
\begin{align*}
    \frac{\partial u}{\partial x_1} & = \sqrt{f(u,\bm{x},\theta)-(\frac{\partial u}{\partial x_2})^2},\\
    \frac{\partial u}{\partial x_1} & = -\sqrt{f(u,\bm{x},\theta)-(\frac{\partial u}{\partial x_2})^2},
\end{align*}
both of whose solutions are solutions of (\ref{counter_eg}). 
Thus, without additional information (eg. $u$ is a function increasing with $x_1$), there is no unique form satisfying (\ref{eq:pde_operator_nonlinear}) that is equivalent to (\ref{counter_eg}). 
\end{example}

\subsection{Incorporating Boundary/Initial Conditions}\label{subsec:boundary_condition}
From both theoretical and numerical points of view, initial and boundary conditions (IBCs) are vital to ensure the uniqueness of a PDE solution \cite{busse2017boundary}. 
Meanwhile, in many practical problems, IBCs are available from physical knowledge or experimental setup without extra experimental cost. 
One example is found in Burger's equation \eqref{eq:burgersequation} (see Section \ref{subsec:eg_burger_eq} for more details), where the PDE is used to model the speed of the fluid, the value of which should be $0$ (Dirichlet boundary) on the adiabatic boundary. 
The information provided by IBCs can help to improve the performance of the GP emulator if they are taken into account during the GP model construction \cite{tan2018gaussian, ding2019bdrygp}. 
Thus, it is important to gain information from IBCs especially when the data are sparse.

In this section, we propose an idea to incorporate the IBCs of PDE into the PIGP method.
In particular, we divide the IBCs into two categories: Dirichlet IBCs and Non-Dirichlet IBCs.
Dirichlet IBCs are given by the known value of PDE solution on specific boundary regions, i.e.,
\begin{equation}\label{eq:dirichlet_boundary}
	u(\bm{x})=b_1(\bm{x}), \bm{x}\in \Gamma_1.
\end{equation}
where $\Gamma_1\subset\partial \Omega$ is the (subset of) boundary region of $\Omega$. 
In many applications, the initial conditions (defined by the known value of PDE solution at time $t = 0$ for transient PDE) of PDEs have the same form as the Dirichlet boundary condition. 
Non-Dirichlet boundary conditions can be represented by a differential operator
\begin{equation}\label{eq:other_boundary}
	\mathcal{B}_{\bm{x}} u(\bm{x})=b_2(\bm{x}, u(\bm{x}), \bm{\theta}), \bm{x}\in \Gamma_2,
\end{equation}
where $\mathcal{B}$ is a differential operator with order $b>0$ and has the similar form with $\mathcal{L}_{\bm{x}}$ that we defined in previous sections, and where $\Gamma_2\subset\partial \Omega$. 
One example of a non-Dirichlet boundary condition is the Neumann boundary condition, i.e., $\nabla u(\bm{x})\cdot \bm{\omega} = b(\bm{x})$, where $\bm{\omega}$ is the outer normal vector. 
Non-Dirichlet boundary conditions can be treated as a part of (\ref{eq:general_pde}).
Thus, we define the comprehensive PDE operator 
\begin{equation*}
	\mathcal{A}(u,\bm{x},\bm{\theta})=\left\{ \begin{aligned}
	&\mathcal{A}(u,\bm{x},\bm{\theta}), & \bm{x}\in \Omega \backslash \Gamma_2  \\
	&\mathcal{B}_{\bm{x}}u(\bm{x}),  & \bm{x}\in \Gamma_2
\end{aligned} \right..
\end{equation*}
The comprehensive PDE operator is obtained by replacing the original PDE operator with the boundary operator $\mathcal{B}_{\bm{x}}u(\bm{x})$ on the boundary region. 
Similarly, we define the comprehensive zeroth-order term
\begin{equation*}
	f(\bm{x}, u(\bm{x}), \bm{\theta})=\left\{ 
	\begin{aligned}
		&f(\bm{x}, u(\bm{x}), \bm{\theta}), & \bm{x}\in \Omega \backslash \Gamma_2 \\
		&b_2(\bm{x}, u(\bm{x}), \bm{\theta}),  & \bm{x}\in \Gamma_2
	\end{aligned} \right..
\end{equation*}
Let $\bm{I}_1\subset \Gamma_1$ and $\bm{I}_2\subset \Gamma_2$ denote the discretization subset of $\Gamma_1$ and $\Gamma_2$, respectively. 
For non-Dirichlet boundary conditions, it is natural to incorporate the boundary information by replacing $\bm{I}$ with $\bm{I} \cup \bm{I}_2$ and replacing $\mathcal{L}_{\bm{x}}$ and $f$ with their comprehensive forms. 
For the Dirichlet boundary condition, we assume that the PDE solution can be observed directly on $\Gamma_1$ without random error. 
Thus, the posterior is now modified as follows,
\begin{align}\label{eq:general_posterior_bound_modify}
	&p\left(\sigma_e^2,\bm{\theta}, u(\bm{I})|W_{\bm{I}}=0, Y(\bm{\tau})=y(\bm{\tau}), U(\bm{I}_1) = b_1(\bm{I}_1)\right)\notag \\
	\propto &\pi(\sigma_e^2)\times \pi_{\bm{\Theta}}\left(\bm{\theta}\right)\times P\left(U(\bm{I})=u(\bm{I})\right)\notag\\
	&\times P\left(Y(\bm{\tau})=y(\bm{\tau}), U(\bm{I}_1)=b_1(\bm{I}_1)|U(\bm{I})=u(\bm{I}),\pi_{\bm{\Theta}}\right)\notag\\
	&\times P\left(W_{\bm{I}}=0|U(\bm{I}_1)=b_1(\bm{I}_1), U(\bm{I})=u(\bm{I}), \bm{\Theta}=\bm{\theta}\right).\notag\\
	\propto  & \frac{1}{\sigma_e^2}\times \exp \Big\{-\frac{1}{2}\big[\left\| u(\bm{I})-\mu(\bm{I})\right\|_{C^{-1}}\notag\\ 
	&+n\log(\sigma_e^2)+\left\| y(\bm{\tau})-u(\bm{\tau})\right\|_{\sigma_e^{-2}}+ \left\|b_1(\bm{I}_1)-\mu(\bm{I}_1) - m_1(u(\bm{I})-\mu(\bm{I}))\right\|_{C_b^{-1}}\notag\\ 
	&+\left\|f(\bm{I},u(\bm{I}),\bm{\theta})-\mathcal{L}_{\bm{x}}\mu(\bm{I})-m_b\{u(\bm{I}\cup\bm{I}_1))-\mu(\bm{I}\cup\bm{I}_1))\} \right\|_{K_b^{-1}}\big]  \Big\},
\end{align}
where $C_b = \mathcal{K}(\bm{I}_1,\bm{I}_1)- \mathcal{K}(\bm{I}_1,\bm{I})\mathcal{K}(\bm{I},\bm{I})^{-1} \mathcal{K}(\bm{I},\bm{I_1})$, $K_b = \mathcal{LKL}(\bm{I},\bm{I})-\mathcal{LK}(\bm{I},\bm{I}\cup\bm{I}_1)\mathcal{K}(\bm{I}\cup\bm{I}_1,\bm{I}\cup\bm{I}_1)^{-1}\mathcal{KL}(\bm{I}\cup\bm{I}_1,\bm{I})$ and $m_b = \mathcal{LK}(\bm{I},\bm{I}\cup\bm{I}_1))\mathcal{K}(\bm{I}\cup\bm{I}_1),\bm{I}\cup\bm{I}_1))^{-1}$, $m_1 = \mathcal{K}(\bm{I}_1,\bm{I})\mathcal{K}(\bm{I},\bm{I})^{-1}$.
To balance $\bm{I}$, $\bm{I}_1$ and $\bm{\tau}$, we propose applying the tempering idea on Dirichlet IBCs, a detailed discussion can be found in Section \ref{subsec:tempering}.
\begin{remark}\label{rmk:ibcs}
Unlike most numerical algorithms for solving PDE, such as finite element methods (FEM) and finite difference methods (FDM), which \textit{must} specify IBCs to ensure the uniqueness of the numerical solution, the PIGP method does not make IBCs mandatory. 
Thus, when IBCs are unavailable, parameter inference methods involving numerically solving PDE using FEM or FDM will not be directly applicable.
To apply these methods, one has to adopt a numerical solution that does not require IBCs, such as collocation methods. 
Otherwise, additional assumptions about the IBCs, such as their parametric forms, must be made to mitigate this problem.
On the other hand, PIGP naturally accommodates PDE parameter inference without IBCs.
This is demonstrated in Section \ref{subsec:toyeg1}.
\end{remark}

\section{Algorithms}\label{sec:algorithm}
We summarize the details for implementing PIGP in this section. 
We first give an overview for PIGP parameter inference and uncertainty quantification. 
Then, we introduce the detailed, step-by-step procedures in the following subsections. 
Our method is implemented in Python and all codes for reproducing
the results shown in this paper are posted on GitHub \footnote{\url{https://github.com/ustclzh/PIGPI}}.
\subsection{Algorithm Summary}\label{subsec:algorithm}
The workflow for PIGP method is summarized as Algorithm \ref{alg:general_alg} in this section. 
\begin{algorithm}[H]\label{alg:general_alg}
\caption{PDE Parameter Inference via PIGP Procedure}
\begin{algorithmic}[1]
\STATE{Given input data $\bm{y}_{\bm{\tau}}$, train GP hyper-parameters based on $\bm{y}_{\bm{\tau}}$.}
\STATE{If PDE operator is linear and independent to $\bm{\theta}$, skip Step 3 and go to Step 4; otherwise go to Step 3.} 
\STATE{Obtain augmented PDE according to Section \ref{subsec:PIGPI2}, initialize GP prior model for all components in augmented PDE.}
\STATE{Construct $\bm{I}$ using method proposed in Section \ref{subsec:determine_I}. If no IBCs are available, skip Step 4 and go to Step 5 with the posterior density (\ref{eq:general_posterior}); otherwise go to Step 6.}
\STATE{Construct $\bm{I}_1$, $\bm{I}_2$ for IBCs and then obtain the posterior density (\ref{eq:general_posterior_bound_modify}) as discussed in Section \ref{subsec:boundary_condition}.} 
\STATE{Optimize the posterior density obtained in previous steps, obtain MAP estimation of $(\sigma_e^2, u(\bm{I}), \bm{\theta})$.}
\STATE{Draw posterior sample for $(\sigma_e^2, u(\bm{I}), \bm{\theta})$ using HMC algorithm or calculate the normal approximation for $(\sigma_e^2, u(\bm{I}), \bm{\theta})$.}
\end{algorithmic}
\end{algorithm}
Steps 2 and 3 have been discussed in Section \ref{sec:method}. 
In the remainder of this Section, we will discuss Steps 4 (Section \ref{subsec:determine_I}), 5 (Section \ref{subsec:map_estimation}), and 6 (Section \ref{subsec:uq_posterior}) in detail.

\subsection{Constructing the Discretization Set \texorpdfstring{$\bm{I}$}{I}}\label{subsec:determine_I}
As stated in Section \ref{subsec:PIGPI1}, the main purpose of introducing $\bm{I}$ is to approximate $W$ using $W_{\bm{I}}$. 
Intuitively, $\bm{I}$ should be dense and space-filling within $\Omega$. 
However, the computational cost will increase with large $n_{\bm{I}} = |\bm{I}|$.
To address this trade-off, we propose a sequential strategy to generate $\bm{I}$ that contains $\bm{\tau}$. 
In practical applications, it is common for the observation data to be collected ahead of data analysis. 
Thus, we assume $\bm{\tau}$ is known and fixed. 
\begin{algorithm}[H]\label{alg:generating_I}
\caption{Generating Discretization Set $\bm{I}$ Encompassing $\bm{\tau}$}
\begin{algorithmic}[1]
\STATE{Input: Observation set $\bm{\tau}$, size for discretization set $n_{\bm{I}} = |\bm{I}|$.}
\STATE{Start with $\bm{I} = \bm{\tau}$. Choose a large candidate point set $\mathcal{C}\subset \Omega$ using Latin hypercube design (LHD) algorithms.
Our empirical investigation suggests setting the size for the candidate point set as $|\mathcal{C}| = 20n_{\bm{I}}$}.
\STATE{Find point $\bm{x}^*\in \mathcal{C}$, such that $\bm{x}^* = \text{argmax}_{\bm{x}\in \mathcal{C}} d(\bm{x},\bm{I})$, where $d(\bm{x},\bm{I}) = \min_{\bm{x}'\in\bm{I}}\|\bm{x}-\bm{x}'\|$ denote the distance of $\bm{x}$ and $\bm{I}$. Let $\bm{I} = \bm{I}\cup \bm{x}^*$ 
}
\STATE{Repeat Step 3 until stopping criterion (Size of $\bm{I}$ equals to $n_{\bm{I}}$.) is reached.}
\end{algorithmic}
\end{algorithm}
The key steps for this algorithm are Steps 2 and 3. 
In Step 2, we generate a sufficiently large candidate set for selecting the design. 
In Step 3, we sequentially select design points that maximize the separation distance of the target design \cite{santner2019design}.
There are two important aspects of using Step 2 to generate a candidate set, instead of directly selecting the points $\bm{x}'$ in $\Omega$ that maximize $d(\bm{x}',\bm{I})$. 
First, when geometry of $\Omega$ is very complicated, optimization of $\text{argmax}_{\bm{x}\in \Omega} d(\bm{x},\bm{I})$ can be challenging \cite{tan2013minimax}. 
Second, by selecting a candidate set from an LHD, we can ensure that the resulting $\bm{I}$ has a non-overlapped one-dimensional projection, which is important for the GP model when the variations of function at each dimension are unknown.

\subsection{Dimension Reduction for \texorpdfstring{$U(\bm{I})$}{U(I)}} \label{subsec:dimension_reduction}
The posterior density (\ref{eq:general_posterior}) has $k\times n_{\bm{I}}+d+k_0|\bm{\tau}|$ ($l\times n_{\bm{I}}+d+k_0|\bm{\tau}|$ for augmented method) variables, where $k$ and $k_0$ are numbers of components and number of observed components in original PDE  
, $d$ is the dimensionality of $\bm{\theta}$. 
To ensure the accuracy of approximation $W_{\bm{I}}$, we intend to select a large discretization set $\bm{I}$ within the computation budget.
Thus, this section seeks to further reduce the computational cost in optimizing or drawing posterior samples from the density (\ref{eq:general_posterior_bound_modify}). 
We propose a dimension reduction method based on the Karhunen-Loève (KL) expansion for the GP model $U(\bm{x})$. 
Recalling that $U(\bm{x})\sim GP(\mu(\bm{x}),\mathcal{K}(\bm{x},\bm{x}'))$, the KL expansion of $U(\bm{x})$ is given by $U(\bm{x})=\mu(\bm{x})+\sum_{i=1}^\infty Z_i \sqrt{\lambda_i}\psi_i(\bm{x})$, where $\lambda_i, \psi_i(\bm{x})$, $i=1,2,\dots$ represent $\mathcal{K}(\bm{x},\bm{x}')$'s eigenvalues and their corresponding eigenfunctions, $\lambda_1\geq \lambda_2\geq\cdots >0$. 
Due to Mercer's theorem, $\sum_{i=1}^\infty\lambda_i=\int_{\bm{x}\in\Omega} \mathcal{K}(\bm{x},\bm{x})d\bm{x}< \infty$, $\lambda_i\to 0$ as $i\to \infty$. 
Thus, we can choose $M\in \mathbb{N}$ such that $\lambda_i$ for $i>M$ are negligible, then the GP $U(\bm{x})$ is approximately represented by $ U(\bm{x})\approx \mu(\bm{x})+ \sum_{i=1}^M Z_i \sqrt{\lambda_i}\psi_i(\bm{x}).$
In practice, since only a discretized version $U(\bm{x})$ is considered, the KL expansion is reduced to principal component decomposition. 
In particular, we can approximate $U(\bm{I})$ by $ U(\bm{I})\approx \mu(\bm{I}) + \sum_{i=1}^M Z_i \sqrt{\lambda_i}\bm{\varphi}_i$, where $\lambda_i, \bm{\varphi}_i$, $i=1,2,\dots,n_{\bm{I}}$ are $C$'s eigenvalues and corresponding eigenvectors.
$M$ can be chosen such that $\sum_{i=1}^M\lambda_i/\sum_{i=1}^{n_{\bm{I}}}\lambda_i\geq 99.99\%$. 
There are two main advantages of applying KL expansion.
First, with the KL expansion, the computation cost of the quadratic form in Eq.\eqref{eq:general_posterior_bound_modify} is reduced to $\mathcal{O}(Mn_I)$, compared to the original cost of $\mathcal{O}(n_I^2)$.
Appendix \ref{spsubsec:kl_expansion_remark} has computational details for the posterior density of KL coefficients $\bm{Z}$ and parameters $\bm{\theta}$. 
Second, the KL expansion transfers the correlated random variable $U(\bm{I})$ to a set of independent random variables $\bm{Z}$, which is helpful for posterior sampling/optimization. 

\subsection{Tempering}\label{subsec:tempering}
As mentioned in Section \ref{subsec:determine_I} and \ref{subsec:boundary_condition}, it is necessary to choose a large discretization set $\bm{I}$ that guarantees $W_{\bm{I}}$ can accurately approximate $W$.
However, as we increase the size of $\bm{I}$, the contribution of GP prior to the posterior will increase while the contribution of the likelihood (i.e., observation data) keeps invariant.
Moreover, when Dirichlet boundary conditions are available, there is also an issue of imbalance between the contributions of boundary conditions and likelihood. 
In order to balance the contribution from different components, in this section we propose the tempering idea. \\
\textbf{Prior Tempering:}
We propose using prior tempering to balance the contribution from the GP prior and that from the observed data. 
First, we ignore the initial boundary conditions herein and we will incorporate IBCs in the next paragraph. 
We can rewrite posterior density (\ref{eq:general_posterior}) in the form 
\begin{align*}
    &p_{\bm{\Theta}, U(\bm{I})| W_{\bm{I}}, Y(\bm{\tau})}\left(\bm{\theta}, u(\bm{I})|W_{\bm{I}}=0, Y(\bm{\tau})=y(\bm{\tau})\right)\\
    \propto & p_{\bm{\Theta}, U(\bm{I})| W_{\bm{I}}}\left(\bm{\theta}, u(\bm{I})|W_{\bm{I}}=0\right)
    p_{Y(\bm{\tau})| U(\bm{I})}\left(y(\bm{\tau})|u(\bm{I})\right).
\end{align*}
As the size of discretization set $|I|$ increases, the prior part $p_{\bm{\Theta}, U(\bm{I})| W_{\bm{I}}}\left(\bm{\theta}, u(\bm{I})|W_{\bm{I}}=0 \right)$ grows, while the likelihood part $p_{Y(\bm{\tau})| U(\bm{I})}\left(y(\bm{\tau})|u(\bm{I})\right)$ stays unchanged. 
Thus, to balance the influence of the prior, we introduce a tempering hyper-parameter $\beta$ and adjust the posterior accordingly
\begin{align*}
    &p^{(\beta)}_{\bm{\Theta}, U(\bm{I})| W_{\bm{I}}, Y(\bm{\tau})}\left(\bm{\theta}, u(\bm{I})|W_{\bm{I}}=0, Y(\bm{\tau})=y(\bm{\tau})\right)\\
    \propto &\pi(\bm{\theta}) \times p_{U(\bm{I})| W_{\bm{I}}}\left(u(\bm{I})|W_{\bm{I}}=0 \right)^{1/\beta}
    p_{Y(\bm{\tau})| U(\bm{I})}\left(y(\bm{\tau})|u(\bm{I})\right).
\end{align*}
A useful setting that we recommend is $\beta=\frac{l|I|}{n}$, where $l$ is the number of system components (for augmented PDE, it is the number of components in augmented PDE), $|\bm{I}|$ is the number of discretization time points, and $n$ is the total number of observation data.
This setting aims to balance the likelihood contribution from the observations with the total number of discretization points.\\
\textbf{IBC Tempering:}
After incorporating the IBCs, the posterior is modified as (\ref{eq:general_posterior_bound_modify}).
The choice of $\bm{I}_1$ can significantly influence the posterior density from two aspects.
First, when $\bm{I}_1$ is a large set while $\bm{\tau}$ is relatively small, i.e., the observation is sparse, the size of the observation set and the boundary set is imbalanced. 
Second, changing the size of $\bm{I}_1$ will significantly change the importance of Dirichlet IBCs in posterior density.
To mitigate these issues, we use a similar tempering idea to balance the influence of the Dirichlet IBCs. 
In particular, we modify the next-to-last line of \eqref{eq:general_posterior_bound_modify} by multiplying  
$\frac{1}{n_1}$ on it, i.e., the next-to-last line of \eqref{eq:general_posterior_bound_modify} becomes $\log(2\pi) + \frac{1}{n_1}\log|C_b| + \frac{1}{n_1}\left\|b_1(\bm{I}_1)-\mu(\bm{I}_1) - \mathcal{K}(\bm{I}_1,\bm{I})\mathcal{K}(\bm{I},\bm{I})^{-1}(u(\bm{I})-\mu(\bm{I}))\right\|_{C_b^{-1}}$.

In summary, the tempering, when Dirichlet IBCs are available, is obtained by replacing the posterior (\ref{eq:general_posterior_bound_modify}) with its modified form 
\begin{align*}
	&p^{(\beta)}_{\sigma_e^2, \bm{\Theta}, U(\bm{I})| W_{\bm{I}}, Y(\bm{\tau}), U(\bm{I}_1)}\left(\sigma_e^2,\bm{\theta}, u(\bm{I})|W_{\bm{I}}=0, Y(\bm{\tau})=y(\bm{\tau}), U(\bm{I}_1) = b(\bm{I}_1)\right)\notag \\
	\propto & P\left(Y(\bm{\tau})=y(\bm{\tau})|U(\bm{I})=u(\bm{I}),\bm{\Theta}\right)\notag\\
	&\times \frac{1}{\sigma_e^2} \times \left[P\left(U(\bm{I})=u(\bm{I})\right)\times P\left(W_{\bm{I}}=0|U(\bm{I}_1)=b(\bm{I}_1), U(\bm{I})=u(\bm{I}), \bm{\Theta}=\bm{\theta}\right)\right]^{\frac{1}{\beta}}\notag\\
	&\times P\left(U(\bm{I}_1)=b(\bm{I}_1)|U(\bm{I})=u(\bm{I})\right)^{\frac{1}{n_1}}\notag
\end{align*}
The tempering balances the posterior contributions from the likelihood, the number of discretization points $|\bm{I}|$, and also the number of  discretization points on the boundary $|\bm{I}_1|$.
\subsection{Maximum a Posteriori Estimation}\label{subsec:map_estimation}
After specifying the unnormalized posterior density of $\bm{\theta}, U(\bm{I})$ and $\sigma_e^2$ as in (\ref{eq:general_posterior}), the maximum a posteriori (MAP) estimation is immediately obtained by maximizing the logarithm of unnormalized posterior density, i.e., the MAP is obtained by
\begin{align}\label{eq:map_opt_sig_e}
	(\hat{\bm{\theta}}, \hat{u}(\bm{I}), \hat{\sigma}_e^2)=&\text{argmax}_{\bm{\theta},u(\bm{I}), {\sigma_e^2} } 	\Big\{- \underbrace{(n+2)\log(\sigma_e^2) - \left\| u(\bm{\tau})-y(\bm{\tau})\right\|_{\sigma_e^{-2}}}_{\text{(log)Likelihood}}\\
	&-\underbrace{\frac{1}{\beta}(\left\| u(\bm{I})-\mu\right\|_{C^{-1}} +\left\|f(\bm{I},u(\bm{I}),\bm{\theta})-\mathcal{L}_{\bm{x}}\mu-m\{u(\bm{I})-\mu\} \right\|_{K^{-1}})}_{\text{Tempered (log)Prior}}\notag\\
	&-\underbrace{\frac{1}{n_1}\left\|b_1(\bm{I}_1)-\mu(\bm{I}_1) - \mathcal{K}(\bm{I}_1,\bm{I})\mathcal{K}(\bm{I},\bm{I})^{-1}(u(\bm{I})-\mu(\bm{I}))\right\|_{C_b^{-1}}}_{\text{Tempered (log)IBCs}}\Big\}.
\end{align}
Note that in some practical applications, the precise information on the measurement devices is available, resulting in a known random error variance $\sigma_e^2$. 
In this case, the MAP estimation for $\bm{\theta}, u(\bm{I})$ is obtained by minimizing (\ref{eq:map_opt_sig_e}) with respect to $\bm{\theta}, u(\bm{I})$ while assuming $\sigma_e^2$ is known. 
Appendix \ref{spsubsec:map_opt_alg} has a detailed discussion about the choice of optimization algorithm for solving (\ref{eq:map_opt_sig_e}).\par
\begin{remark}
Due to the fact that $\bm{I}$ is a discretization set, which is a finite point set, the estimation of PDE solution $u(\bm{x})$ at point $\bm{x}\notin \bm{I}$ is obtained by the posterior mean of the conditional GP, $U(\bm{x})|U(\bm{I}) = \hat{u}(\bm{I})$. 
Appendix \ref{spsubsec:predict_GP} has details for obtaining the estimation of $U(\bm{x})$ for all $\bm{x}\in \Omega$.
\end{remark}
\subsection{Simultaneous Uncertainty Quantification for Parameters and PDE Solution}\label{subsec:uq_posterior}
To quantify the uncertainty of unknown parameters, we recommend two methods to approximate posterior density. 
One is to use multivariate normal approximation, and the other is to sample from the posterior density via a Hamiltonian Monte Carlo (HMC).\\
\textbf{Normal Approximation for Posterior Density \cite{gelman2013bayesian}:} Let $(\hat{\bm{\theta}}, \hat{u}(\bm{I}),\hat{\sigma}_e^2)$ denote the MAP of $(\bm{\theta}, U(\bm{I}),\sigma_e^2)$, we can calculate the hessian matrix $H$ of minus log-posterior density (Unnormalized) at $(\hat{\bm{\theta}}, \hat{u}(\bm{I}),\hat{\sigma}_e^2)$. 
The posterior is thus approximated by $N((\hat{\bm{\theta}}, \hat{u}(\bm{I}), \hat{\sigma}_e^2), H^{-1})$.
The component-wise $1-\alpha$ credible interval for $\theta_i$ can be obtained as $[\hat{\theta}_i-\Phi^{-1}(\alpha/2)h_i,\hat{\theta}_i+\Phi^{-1}(1-\alpha/2)h_i]$, where $\hat{\theta}_i$ is the $i-$th element of $\hat{\bm{\theta}}$, $\Phi^{-1}(\cdot)$ is the inverse of distribution function for the standard normal random variable, $h_i^2$ is the $i-$th diagonal element of $H^{-1}$.
Credible intervals for $\sigma_e^2$ and $u(\bm{I})$ are obtained similarly. \\
\textbf{Hamiltonian Monte Carlo\cite{neal2011mcmc}} can be employed to draw a posterior sample of $U(\bm{I})$, $\bm{\theta}$ and $\sigma_e^2$. 
To speed up the burn-in process, we suggest initializing HMC using the MAP estimation of $\bm{\theta}, U(\bm{I})$ and $\sigma_e^2$ obtained from Section \ref{subsec:map_estimation}. (Appendix \ref{spsubsec:uq_posterior_remark} has details for implementation of HMC and tempering of the posterior density.).\par
One of the most commonly used measures of uncertainty in Bayesian inference is the credible interval.
When an HMC sample is obtained, one can easily approximate the point-wise credible interval for parameters using the sample quantiles of their posterior sample.

\section{Numerical Illustrations}\label{sec:simulation_res}F
In this section, we employ four examples to illustrate the accuracy and efficiency of PIGP.
Before showing the simulation results, we first briefly introduce the benchmark methods and evaluation metrics.\\ 
\textbf{Benchmark Methods}: For demonstration purposes, we choose several methods as benchmarks. 
We illustrate the performance of the proposed method by comparing it with these benchmark methods.
\begin{itemize}
    \item The first benchmark method is maximum likelihood estimation (MLE). 
    When we ignore the computational complexity of solving a PDE, the MLE can be obtained by maximizing the log-likelihood $l(\bm{\theta},\sigma_e^2) = \frac{1}{2}\log(\sigma_e^2) -\frac{1}{2\sigma_e^2}\sum_{i=1}^n (y_i-u(\bm{x}_i,\bm{\theta}))^2$ using state-of-the-art optimization algorithms, here we use the notation $u(\bm{x}_i,\bm{\theta})$ to denote the PDE solution at input variable $\bm{x}_i$ given parameter $\bm{\theta}$, obtained by numerical solution of the PDE. 
    It is easy to see that the MLE of $\bm{\theta}$, denoted by $\bm{\hat{\theta}}$, can be obtained by minimizing $L(\bm{\theta}) = \sum_{i=1}^n (y_i-u(\bm{x}_i,\bm{\theta}))^2$ and MLE of $\sigma_e^2$ is $\frac{1}{n}\sum_{i=1}^n (y_i-u(\bm{x}_i,\bm{\hat{\theta}}))^2$.
    \item The second benchmark method is the two-stage method (TSM) \cite{rai2019gaussian}. 
    The two-stage method for PDE parameter inference is constructed as follows. 
    In the first stage, the PDE solution $u(\bm{x})$ and its partial derivatives are estimated by the Kriging model and its derivatives. 
    In the second stage, we plug in the Kriging model and its derivatives in PDE, the parameter $\bm{\theta}$ are estimated through the least squared method that minimizes the error $\sum_{i=1}^{n_I}(\mathcal{A}(u(\bm{x}_i),\bm{x}_i,\bm{\theta})-f(\bm{x}_i,u(\bm{x}_i),\theta))^2$. 
    This procedure is shown to be slightly better than the two-stage method employed in \cite{xun2013parameter}, as shown in Section \ref{subsec:eglidar}. 
    \item The third benchmark method is the automated PDE identification (API) method proposed in \cite{liu2021automated}.
    The API method is proposed for PDE identification, which can also be used for PDE parameter inference.
    \item The fourth method is the Bayesian optimization method (BOM). 
    The method is based on the pioneer work proposed in \cite{jones1998efficient}, which is also called the expected improvement method.
    We construct BOM to optimize $\log(L(\bm{\theta}))$.
    The BOM is well known for its efficiency and accuracy. 
    For BOM, we set the budget of numerical simulations to $20p$, where $p$ is the dimension of parameter space.
    The budget that we set is larger than suggested in \cite{loeppky2009choosing}, which is helpful to improve the performance of BOM.
    \item The fifth benchmark method is proposed by \cite{xun2013parameter}.
    This benchmark method consists of two methods, called the Bayesian method (Xun-BM) and the parameter cascading method (Xun-PC). 
    The performance of these two methods is directly taken from \cite{xun2013parameter}, based on the LIDAR equation, and we will compare these results with PIGP in Section \ref{subsec:eglidar}.
\end{itemize}
\textbf{Evaluation metrics}: We consider multiple metrics with the main purpose listed below among multi-run of simulation:
\begin{itemize}
    \item Root mean squares errors (RMSEs), the biases, and the standard deviation (SD) of the estimated parameters: for the accuracy of parameter estimates.
    \item RMSEs of estimation for $u(\bm{x})$ on a test set $\bm{I}_p$: the quality of PDE solution recovery.
    \item Coverage rates (CR) of $95\%$ credible intervals of posterior: the performance of uncertainty quantification of Bayesian posterior.
    \item The computational time of PIGP with augmented PDE method, compared to both the PIGP without augmented PDE method (Section \ref{subsec:eglidar}) and other benchmark methods: the efficiency of employing augmented PDE.\footnote{For fairness in computational time comparison, we record the computational time of all the codes on a desktop with Intel Xeon CPU E5-1650 v3 (3.50GHz) and 32 GB RAM.
    For all other performance comparisons, the codes are deployed on a high-performance computer cluster.}
\end{itemize}
Now, we briefly introduce the examples used in this section and highlight the rationale behind each selection.\\
\textbf{Highlights in Each Example:}
In the first example, we introduce a toy example simplified from the disease transmission model with age structure to illustrate the performance of PIGP when IBCs are unavailable. 
The second example of the LIDAR equation shows the efficiency of PIGP with augmentation by comparing it with the PIGP with the original PDE. 
Furthermore, we illustrate the accuracy of PIGP by comparing it with the methods proposed in \cite{xun2013parameter} and the two-stage method. 
The third example is based on Burger's equation, which is employed to illustrate the performance of PIGP on parameter inference for nonlinear PDEs. 
In the last example, a multi-dimensional PDE called the diffusion-Brusselator equation is introduced to demonstrate the performance of PIGP on parameter inference and PDE solution recovery when the observation data is censored. \\
\textbf{Synthetic Observation Data Generation Mechanism:} 
In order to assess the performance under different scenarios, we need to generate a multi-set of synthetic data. 
Each set of synthetic data is generated using the following mechanism:
First, we determine a true value for $\bm{\theta}$, denoted by $\bm{\theta}_0$. 
Then we generate $\bm{\tau}= \{\bm{x}_1,\dots,\bm{x}_n\}\subset \Omega$ from mesh grid or using quasi-Monte Carlo algorithm/space-filling designs. 
The experimental data $\{y_i, i =1,\dots, n\}$ is then generated from PDE solution with given $\bm{\theta}=\bm{\theta}_0$ on $\bm{\tau}$ through adding i.i.d. Gaussian errors, i.e., the data is $y_i = u(\bm{x}_i)+\varepsilon_i, i=1,\dots,n$ where $\varepsilon_i \overset{i.i.d.}{\sim} \mathcal{N}(0, \sigma_e^2)$.
In this paper, all PDEs are solved numerically using MATLAB. 
All MATLAB codes for obtaining numerical PDE solutions are available on the same GitHub.\par
It's noteworthy that not all benchmark methods can be applied to every example presented in this section (with the exception of MLE, which is \textit{theoretically} applicable to all examples. However, the MLE method becomes impractical in scenarios where numerical PDE solving is too time-consuming as it relies on numerous evaluations of the PDE solver). 
We will compare our PIGP with benchmark methods wherever applicable for each individual example.
All codes for reproducing the results shown in this section are also posted on GitHub.
\subsection{Toy Example: Disease Transmission Model with Age Structure} \label{subsec:toyeg1}
We first consider a simple toy example. 
The background of the PDE is the modeling of populations and infectious diseases (see Chap. 13.2 of \cite{brauer2019correction}). 
The model assumes that the infected population can be described by a density function $u(t,a)$ of age $a$ (continuous variable) and time $t$.
In particular, the PDE is given by 
\begin{equation}\label{eq:eq:eg_toy_disease}
	\frac{\partial u(t,a)}{\partial t}+\frac{\partial u(t,a)}{\partial a}=-\mu(a)u(t,a), a\geq 0, t\geq 0
\end{equation}
with boundary and initial conditions $u(t,0)=\int_{0}^{\infty} \beta(a)u(a,t)da=B(t), t\geq 0$,$u(0,a)=u_0(a), a\geq 0$, where $\beta(a)$ is the per-capita age specific fertility rate.
The parameter of interest $\mu(a)\geq 0$ denotes the per-capita age specific death rate.\par
We employ this simple example to demonstrate the performance of the PIGP method in estimating the parameter and PDE solution. 
For simplicity, we consider a parametric form for $f(\bm{x},u(\bm{x}),\bm{\theta}) = -\mu(a)u(\bm{x}) = (\theta_1 +\theta_2 a +\theta_3 a^2)u(a,t)$, where $\bm{x} = (t,a)$.
That is, we assume the unknown $\mu(a)$ function has a parametric form $\mu(a)=-\theta_1 -\theta_2 a -\theta_3 a^2$ . We set the true value for $\bm{\theta}$ at $\bm{\theta}_0 = (1,-2,0)$, so that $\mu(a)=2a-1$, and the true solution of PDE \eqref{eq:eq:eg_toy_disease} is $u(a,t)=\exp\{t-a^2\} $.
The data is observed on a set $\bm{\tau}=\{(t_i,s_i),i=1,\dots,n\}$ by adding random errors to the PDE solution, where $n$ denotes the number of observation data which will be specified for specific simulation cases. 
In this Section, we will examine the performance of alternative methods for $3$ cases, where the size of observation data is set as $n=30,60,120$. 
For each case, we generate $\bm{\tau}$ using maximin LHD \footnote{Implemented by Python scipy package (version 1.8.1)}.
To apply PIGP method, we define $\mathcal{L}_{\bm{x}}=\frac{\partial}{\partial t}+\frac{\partial}{\partial a}$.\par
In this example, we assume the boundary condition is not available, which corresponds to some commonly faced real-world applications where the boundary condition (eg. $B(t)$) and initial condition (eg. $u_0(a)$) are very hard to obtain.
The MLE and BOM rely on the numerical solution of PDE with given parameters, in which IBCs are always required to ensure the uniqueness of the numerical solution.
To apply these two methods, we parameterize the IBCs by assuming the boundary condition has the form $B(t) = \exp(\vartheta_0+\vartheta_1 t)$, and the initial condition has the form $u_0(a) = \exp(\vartheta_2 +\vartheta_3 a + \vartheta_4 a^2)$. 
We adopt these conditions into MLE and BOM.
Note that MLE and BOM estimate the parameters $\bm{\theta}$ and $\bm{\vartheta} = (\vartheta_0,\vartheta_1,\dots,\vartheta_4)$ simultaneously.\par
For each $n$, we repeat the experiment 100 times using each method and calculate the RMSEs of MAP obtained by the PIGP method and parameter estimation obtained by MLE, BOM and TSM. 
The results are shown in Table \ref{tb:compare_toy1}. 
When $n=30$ and $n=60$, PIGP outperforms all other methods in terms of RMSE of parameter estimation, including the MLE. 
Since the MLE and BOM require specified IBCs for the PDE numerical solver, five additional nuisance parameters in IBCs need to be estimated. 
Through IBCs parameterization, MLE and BOM have to estimate much more parameters than PIGP and TSM.
This might reduce the accuracy of parameter estimation by MLE and BOM when the sample size is small.
For the case $n=120$, MLE is slightly better than PIGP, but except that, PIGP still outperforms all other benchmark methods. 
\begin{table}[H]
\caption{The mean±SD of RMSEs of MAPs obtained by PIGP, parameter estimation obtained by TSM, MLE and BOM across 100 simulated datasets. The proposed method is emphasized with boldface font.}\label{tb:compare_toy1}    
\begin{center}
\begin{tabular}{|l|l|l|l|l|l|}
\hline
$n$                 & $n_I$ & \textbf{PIGP}       & TSM         & MLE                          & BOM                        \\ \hline
\multirow{3}{*}{30} & 30    & \textbf{0.030±0.020} & 0.073±0.052 & \multirow{3}{*}{0.040±0.057} & \multirow{3}{*}{3.00±1.03} \\ \cline{2-4}
                    & 60    & \textbf{0.024±0.019} & 0.104±0.068 &                              &                            \\ \cline{2-4}
                    & 120   & \textbf{0.020±0.015} & 0.113±0.064 &                              &                            \\ \hline
\multirow{2}{*}{60} & 60    & \textbf{0.015±0.010} & 0.036±0.024 & \multirow{2}{*}{0.013±0.021} & \multirow{2}{*}{3.13±0.97} \\ \cline{2-4}
                    & 120   & \textbf{0.013±0.009} & 0.057±0.039 &                              &                            \\ \hline
120                 & 120   & \textbf{0.010±0.009} & 0.020±0.012 & 0.006±0.005                  & 3.14±0.93                  \\ \hline
\end{tabular}
\end{center}
\end{table}
We also compare the RMSEs of estimation of PDE solution on a grid $\{0.05,0.15,\dots,0.95\}^2$ obtained by PIGP, TSM, MLE, and BOM. 
For MLE and BOM, the estimation of the PDE solution is obtained by the numerical PDE solution by setting the parameter as the estimated value. 
Note that the estimation of the PDE solution via TSM is equivalent to the Kriging interpolation of the PDE solution using observation data. 
Results listed in Table \ref{tb:est_u_toy} show that PIGP can provide a more accurate estimation for PDE solution $u$ on $\bm{I}$ than Kriging. 
It can be seen that PIGP outperforms TSM and BOM in all cases. 
For MLE, only when $n=120$, MLE can estimate the PDE solution slightly more accurately than PIGP.
It is worth noting that when physical data is very limiting, e.g., $n = 30$ is this example, PIGP outperforms MLE on both estimating parameters $\bm{\theta}$ and PDE solution. 
The reason is that MLE requires estimating $8$ parameters including $5$ nuisance parameter. 
When physical data is not enough, estimation of a large number of parameters is challenging. 

\begin{table}[H]
\caption{The mean±SD of RMSE of PDE solution estimation obtained from the PIGP, TSM, MLE, and BOM. 
All results are rescaled by multiplying $10^{3}$. The proposed method is emphasized with boldface font.}\label{tb:est_u_toy}    
\begin{center}
\begin{tabular}{|r|r|l|l|l|l|}
\hline
\multicolumn{1}{|l|}{$n$} & \multicolumn{1}{l|}{$n_I$} & \textbf{PIGP}     & TSM       & MLE       & BOM     \\ 
\hline
\multirow{3}{*}{30}   & 30    & \textbf{1.39±0.34} & 1.71±0.40 & \multirow{3}{*}{2.56±3.18}  & \multirow{3}{*}{884±193} \\ \cline{2-4} 
 & 60       & \textbf{0.75±0.13} & 1.57±0.36 &          &         \\ \cline{2-4} 
                          & 120                        & \textbf{0.66±0.13} & 1.56±0.35 &           &         \\ \hline
\multirow{2}{*}{60}       & 60                         & \textbf{0.66±0.09} & 0.88±0.11 & \multirow{2}{*}{0.6±0.84}  & \multirow{2}{*}{828±193} \\ \cline{2-4} 
                          & 120                        & \textbf{0.50±0.08} & 0.80±0.11 &           &         \\ \hline
120                       & 120                        & \textbf{0.42±0.05} & 0.57±0.06 & 0.34±0.14 & 815±203 \\ \hline
\end{tabular}
\end{center}
\end{table}
From comparison results in this example, we can conclude that when $n$ is small, i.e., observation is sparse, say $n=30$, choose a large $n_{\bm{I}}$ using the proposed method can significantly reduce the errors of both parameter estimation and PDE solution estimation. 
\subsection{Long-range Infrared Light Detection and Ranging Equation} \label{subsec:eglidar}
In this Section, we revisit the LIDAR equation introduced in Section \ref{subsec:basic_formulation}.
The LIDAR equation is a linear PDE where the unknown parameters are involved in the PDE operator. 
It is employed in \cite{xun2013parameter} to illustrate the performance of their methods on PDE parameter inference. 
In this section, we introduce this example for two purposes: First, by comparing the computational time of PIGP with original PDE and PIGP with augmented PDE, we illustrate the efficiency of the PIGP method with augmentation of PDE. 
Second, through comparing with results proposed in \cite{xun2013parameter} and results of TSM, we illustrate the accuracy of PIGP on estimating parameters. \par
The LIDAR equation is given by \eqref{eq:lidar}, where $t\in[0,20],s\in[0,40]$.
The boundary condition is given by $u(t,0) = u(t,40) = 0$ and initial condition is $u(0,s)=\{1+0.1\times (20-s)^2\}^{-1}$. 
The true value for $\bm{\theta}$ is $\bm{\theta}_0 = (1,0.1,0.1)$. \par
Following the suggestion by \cite{xun2013parameter}, we solve the PDE using MATLAB build-in function \textit{pdpde} on a $40000\times 40000$ grid. 
This configuration requires 1016 seconds to solve the PDE once.
Meanwhile, solving this PDE requires a large amount of RAM (over 12 GB of memory usage is observed). 
For MLE, a global optimization requires at least one thousand evaluations of the PDE solution to guarantee reliable results. 
For BOM, at least $20p = 60$ repeated evaluation is required for each case. 
While in this example, each case is repeated 1000 times. 
Thus, methods that involve repeated evaluation of numerical solutions (i.e., MLE and BOM) are too time-consuming to be practical.\par
\begin{figure}[H]
	\centering \includegraphics[width=15cm]{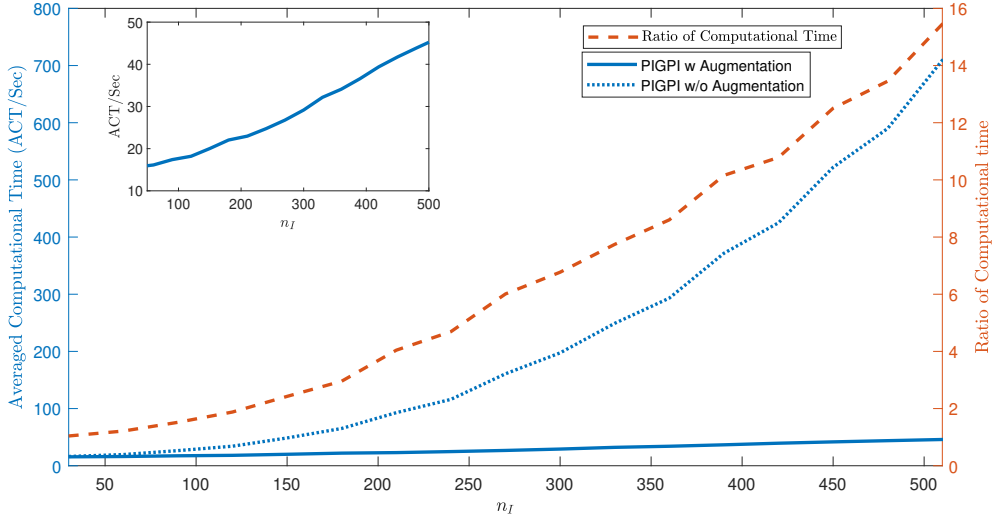}\\
	\centering \caption{The computation time of MAP optimization using PIGP with and without augmented PDE, top left zoom-in figure plots the computation time of MAP optimization using PIGP with augmented PDE for better illustration of its nearly-linear computational time with respect to $n_{\bm{I}}$. (Online version of this figure is colored)}\label{fig:time_comp}
\end{figure}

As shown in Section \ref{subsec:basic_formulation} LIDAR equation is a linear PDE where the PDE parameter depends on the unknown parameters. 
Thus, we employ the PIGP method on augmented PDE.
The augmented PDE for the LIDAR equation that satisfies the ``lowest order of derivative'' principal is given as follows, 
\begin{align*}
    \frac{\partial u_1(t,s)}{\partial s} &= u_2(t,s),\\
    \frac{\partial u_2(t,s)}{\partial s} &= u_3(t,s),\\
    \frac{\partial u_1(t,s)}{\partial t} &= \theta_D u_3(t,s) + \theta_S u_2(t,s) + \theta_A u_1(t,s).
\end{align*}
Due to the linearity of the LIDAR equation, the PIGP based on the original (non-augmented) PDE can also be applied with the caveat in computational cost. 
In this section, we employ the PIGP based on the original PDE as an alternative method. 
We compare PIGP based on augmented PDE and PIGP based on original PDE and show that augmented PDE can significantly reduce the computational time for PDE parameter inference. \\
\textbf{Computational time of PIGP with and without augmentation:} 
As stated in section \ref{subsec:PIGPI2}, when $\mathcal{L}_{\bm{x}}^{\bm{\theta}}$ depends on $\bm{\theta}$, PIGP requires updating covariance matrix $\mathcal{K}$ in each evaluation of posterior density. 
On the other hand, PIGP with augmentation leaves $\mathcal{K}$ unchanged for each evaluation.
We apply the Adam algorithm to optimize the posterior density and set the algorithm to run 2500 iterations.
For specific $n\in\{30,60,\dots, 510\}$, $\bm{\tau}$ is generated using maximin LHD. 
For simplicity we set $\bm{I} = \bm{\tau}$.
The random errors are i.i.d. normal distributed with zero mean and variance $\sigma_e^2 = 0.01^2$. 
Each instance is repeated 100 times and the average computational time of the Adam algorithm is recorded and is shown in Figure \ref{fig:time_comp}. 
It can be seen that the optimization of posterior density via PIGP with augmentation is much faster than the one without. 
Figure \ref{fig:time_comp} also shows that the ratio of computational time between the two methods can significantly increase with the increase of the size of the discretization set. 
Moreover, we compare the RMSE of $\bm{\theta}$ estimation obtained from two alternative methods for $n_{\bm{I}} = 30,60,120,240,480$ in Table \ref{tb:compare_err_aug}.
It can be seen that the accuracy of PIGP parameter estimates with augmentation are better than the accuracy of PIGP parameter estimates without augmentation in most cases.
\begin{table}[H]
\caption{The RMSE±SD of parameter estimation via PIGP and PIGP without augmentation methods. }\label{tb:compare_err_aug}
\begin{tabular}{|l|l|l|l|l|l|}
\hline
$n_I$                                                                             & 30          & 60         & 120        & 240       & 480       \\ \hline
\begin{tabular}[c]{@{}l@{}}PIGP w\\ Augmentation ($\times 10^{-3})$\end{tabular}                  & 41.97±19.72 & 14.87±2.25 & 9.66±3.21  & 6.95±4.26 & 0.55±0.38 \\ \hline
\begin{tabular}[c]{@{}l@{}}PIGP w/o\\ Augmentation ($\times 10^{-3})$\end{tabular} & 91.74±5.47  & 34.84±2.86 & 19.62±2.09 & 6.53±1.06 & 1.02±0.62 \\ \hline
\end{tabular}
\end{table}

In \cite{xun2013parameter}, two methods, called Xun-BM and Xun-PC, are proposed for PDE parameter inference and are illustrated using LIDAR equation. 
In the remainder of this section, we will illustrate the performance of PIGP through comparing it with the BS andXun-PCproposed in \cite{xun2013parameter} and TSM. \\
\textbf{Comparing with Methods Proposed in \cite{xun2013parameter}:}. 
In order to compare with the Xun-BM and Xun-PC methods proposed by Xun et. al. \cite{xun2013parameter}, we conduct the simulation with the same parameter setup with \cite{xun2013parameter}. 
In particular, $\bm{\tau}=\{(i,j),i=1,\dots,20,j=1,\dots,40\}$ is a grid.
Since $n=800$ is large enough, we set $\bm{I} = \bm{\tau}$ in this simulation.
Two cases with different $\sigma_e$ values are evaluated, i.e., $\sigma_e = 0.02,0.05$.
Following the experiment setup in \cite{xun2013parameter}, PIGP and TSM are run 1000 times, respectively. 
The normal approximation for posterior density is employed to calculate the $95\%$ credible interval of $\bm{\theta}$ for PIGP. 
We validate the accuracy of normal approximation by comparing it with HMC results at the end of this section. 
Note that due to the computational burden caused by numerically solving PDE using MATLAB (over 1 hour for evaluating PDE solution once), MLE and BOM are both unavailable. 
As a comparison, PIGP takes less than 2 minutes to solve the PDE parameter inference once. 
The results are shown in Table \ref{tb:eg_lidar_jasa}.
We can see that our method is better than all alternative methods in this comparison. 
\begin{table}[H]
\caption{The mean of biases, standard deviations, component-wise RMSE of the parameter estimates using the proposed PIGP method, Bayesian method (Xun-BM) (from \cite{xun2013parameter}), the parameter cascading method (Xun-PC) (from \cite{xun2013parameter}), and the two-stage method (TSM) in the $1000$ simulated datasets when the data noise has the standard deviation $\sigma_e = 0.02, 0.05$. The proposed method is emphasized with boldface font.}\label{tb:eg_lidar_jasa}    
\begin{center}
\begin{tabular}{|l|l|rrr|rrr|}
\hline  \multicolumn{2}{|l|}{}    & \multicolumn{3}{l|}{$\sigma_e = 0.02$}   & \multicolumn{3}{l|}{$\sigma_e = 0.05$}\\ \hline    \multicolumn{2}{|l|}{}     & \multicolumn{1}{l|}{$\theta_D$}      & \multicolumn{1}{l|}{$\theta_A$}     & \multicolumn{1}{l|}{$\theta_S$} & \multicolumn{1}{l|}{$\theta_D$}      & \multicolumn{1}{l|}{$\theta_A$}     & \multicolumn{1}{l|}{$\theta_S$} \\ \hline
\multirow{4}{*}{\begin{tabular}[c]{@{}l@{}}Bias\\ $\times 10^{-3}$\end{tabular}}   & \textbf{PIGP} & \multicolumn{1}{r|}{\textbf{-14.00}} & \multicolumn{1}{r|}{\textbf{-0.20}} & \textbf{-0.12} & \multicolumn{1}{r|}{\textbf{-27.35}} & \multicolumn{1}{r|}{\textbf{-0.34}} & \textbf{-0.29}  \\ \cline{2-8} 
&Xun-BM & \multicolumn{1}{r|}{-16.50}&\multicolumn{1}{r|}{-0.40} & -0.20   & \multicolumn{1}{r|}{-35.60}  & \multicolumn{1}{r|}{1.00}& 0.60  \\ \cline{2-8} 
&Xun-PC & \multicolumn{1}{r|}{-29.70} & \multicolumn{1}{r|}{-0.10}   & -0.30     & \multicolumn{1}{r|}{-55.90} & \multicolumn{1}{r|}{-0.20} & -0.50 \\ \cline{2-8} 
& TSM  & \multicolumn{1}{r|}{-105.33} & \multicolumn{1}{r|}{-2.69}  & -1.28    & \multicolumn{1}{r|}{-140.12} & \multicolumn{1}{r|}{-4.05} & -2.12  \\ \hline
\multirow{4}{*}{\begin{tabular}[c]{@{}l@{}}SD\\ $\times 10^{-3}$\end{tabular}} & \textbf{PIGP} & \multicolumn{1}{r|}{\textbf{9.37}}   & \multicolumn{1}{r|}{\textbf{1.63}}  & \textbf{0.21}  & \multicolumn{1}{r|}{\textbf{20.31}}  & \multicolumn{1}{r|}{\textbf{3.74}}  & \textbf{0.48}                   \\ \cline{2-8} 
 &Xun-BM & \multicolumn{1}{r|}{9.10} & \multicolumn{1}{r|}{1.60} & 0.20  &\multicolumn{1}{r|}{22.20} & \multicolumn{1}{r|}{3.80} & 0.50  \\ \cline{2-8} 
 &Xun-PC & \multicolumn{1}{r|}{24.90}   & \multicolumn{1}{r|}{3.80}           & 0.50       & \multicolumn{1}{r|}{40.50}    & \multicolumn{1}{r|}{6.20}   & 0.80    \\ \cline{2-8} 
 & TSM  & \multicolumn{1}{r|}{29.42}  & \multicolumn{1}{r|}{3.82}  & 0.52   & \multicolumn{1}{r|}{49.00}  & \multicolumn{1}{r|}{7.29} & 1.03   \\ \hline
\multirow{4}{*}{\begin{tabular}[c]{@{}l@{}}RMSE\\ $\times 10^{-3}$\end{tabular}}  & \textbf{PIGP} & \multicolumn{1}{r|}{\textbf{16.85}}  & \multicolumn{1}{r|}{\textbf{1.64}}  & \textbf{0.24}  & \multicolumn{1}{r|}{\textbf{34.06}}  & \multicolumn{1}{r|}{\textbf{3.75}}  & \textbf{0.56} \\ \cline{2-8} 
 &Xun-BM & \multicolumn{1}{r|}{18.81} & \multicolumn{1}{r|}{1.66} & 0.27 &\multicolumn{1}{r|}{42.00}  & \multicolumn{1}{r|}{3.90}   & 1.00      \\ \cline{2-8} 
&Xun-PC & \multicolumn{1}{r|}{38.96}  & \multicolumn{1}{r|}{3.75} & 0.54   & \multicolumn{1}{r|}{69.10} & \multicolumn{1}{r|}{6.20}  & 2.20    \\ \cline{2-8} 
 & TSM    & \multicolumn{1}{r|}{109.35}   & \multicolumn{1}{r|}{4.67}  & 1.38     &\multicolumn{1}{r|}{148.43}    & \multicolumn{1}{r|}{8.34}   & 2.36 \\ \hline
\multirow{3}{*}{\begin{tabular}[c]{@{}l@{}}CR\\ $\%$\end{tabular}}      & \textbf{PIGP} & \multicolumn{1}{r|}{\textbf{98.6}}   & \multicolumn{1}{r|}{\textbf{100}}   & \textbf{99.2}                   & \multicolumn{1}{r|}{\textbf{79.7}}   & \multicolumn{1}{r|}{\textbf{95.9}}  & \textbf{92.2} \\ \cline{2-8} 
 &Xun-BM & \multicolumn{1}{r|}{93.9}  & \multicolumn{1}{r|}{99.9}   & 98.8    &\multicolumn{1}{r|}{74}  & \multicolumn{1}{r|}{97.8} & 93.5 \\ \cline{2-8} 
 &Xun-PC&\multicolumn{1}{r|}{84.3}&\multicolumn{1}{r|}{96.7}  & 94.9 & \multicolumn{1}{r|}{78.1}  &\multicolumn{1}{r|}{96.5}  & 93.8   \\ \hline
\end{tabular}    
\end{center}
\end{table}

\textbf{Implementing HMC and Validating the Accuracy of Normal Approximation}: We implement the HMC algorithm on the LIDAR example to obtain the posterior sample of $u$ and $\bm{\theta}$.
In particular, we assume that the observation data is on a grid  $\bm{\tau}=\{(i,j),i=1,\dots,20,j=1,\dots,40\}$ and $\bm{I} = \bm{\tau}$. 
The data is polluted by random errors whose distribution is $0$ mean $\sigma_e^2 = 0.02^2$ variance normal distribution.

\begin{figure}[H]
	\centering \includegraphics[width=15cm]{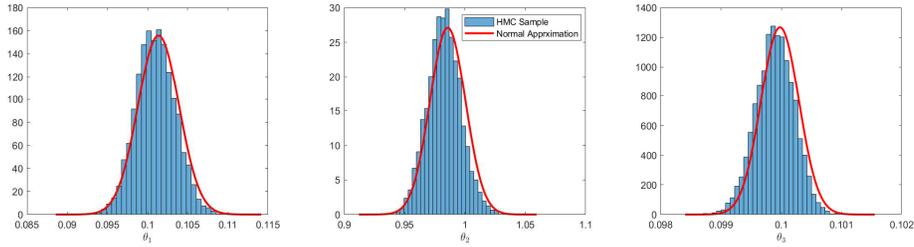}\\
	\centering \caption{LIDAR equation. The posterior sample obtained by HMC, posterior density estimate obtained by normal approximation.}\label{fig:lidar_hmc_normal}
\end{figure}

\begin{figure}[H]
	\centering \includegraphics[width=15cm]{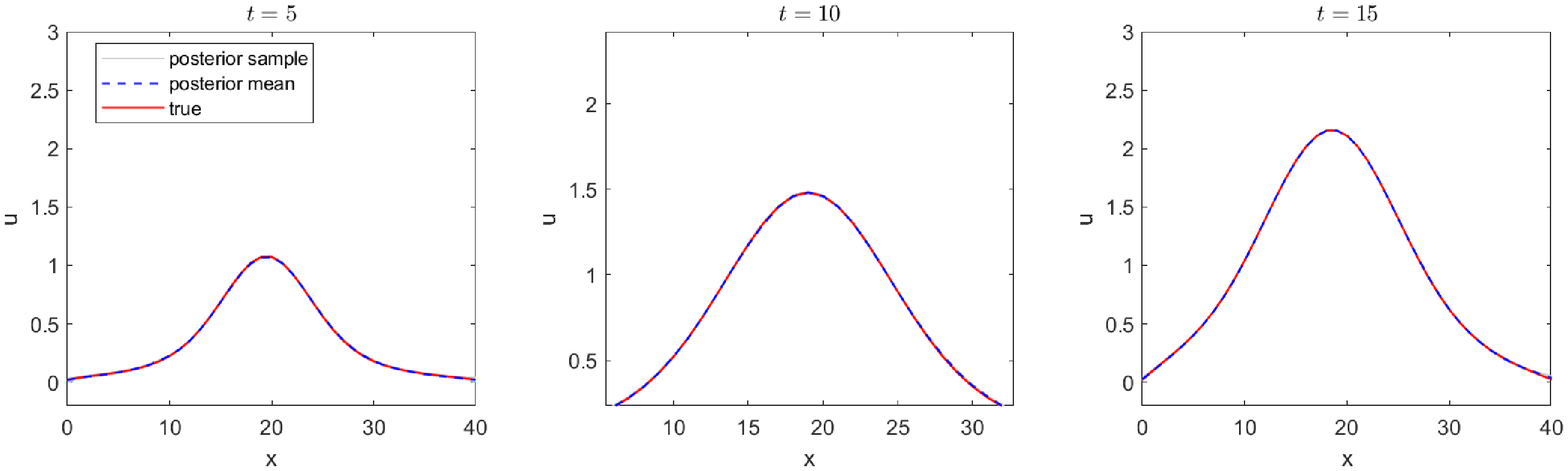}\\
	\centering \caption{LIDAR equation. Posterior sample, posterior mean, and the true value of PDE solutions. The posterior sample and posterior mean are produced by HMC using PIGP, $n=800$.}\label{fig:lidar_hmc}
\end{figure}

\begin{figure}[H]
	\centering \includegraphics[width=15cm]{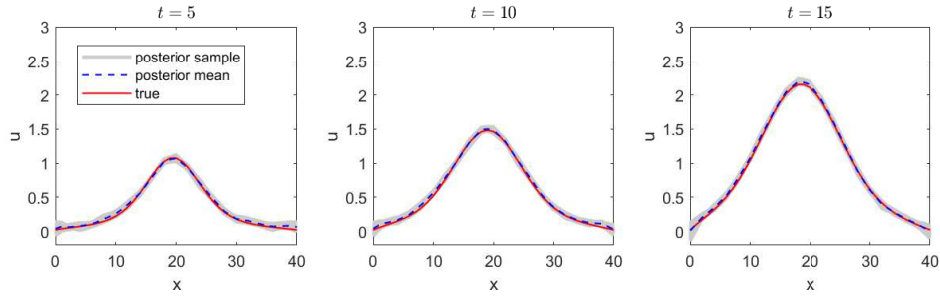}\\
	\centering \caption{LIDAR equation. Posterior sample, posterior mean, and the true value of PDE solutions. The posterior sample and posterior mean are produced by HMC using PIGPP, $n=100$.}\label{fig:lidar_hmc_sparse}
\end{figure}

The main purpose of this simulation is to validate the accuracy of normal approximation by comparing it with the HMC results. 
In this example, the parameters for HMC are tuned such that the acceptance rate is stable between $70\%$ and $90\%$. 
By comparing the posterior density obtained by HMC and normal approximation, as shown in Figure \ref{fig:lidar_hmc_normal}, it can be seen that the normal approximation works well in this example. 
The comparison in this section illustrates the accuracy of the posterior mode and credible interval obtained by normal approximation. 
Figure \ref{fig:lidar_hmc} shows the true value, HMC sample, and posterior mean (estimated from HMC sample) of the PDE solution. 
We plot the profiles of the solution at three time values, i.e., $u(5,s),u(10,s),u(15,s), s\in [0,40]$. 
It can be seen that the HMC sample can cover the true value and the posterior mean is an accurate estimation of the PDE solution.
\blue{Since we set $n=800$ following the setup in \cite{xun2013parameter}, it can be shown in Figure \ref{fig:lidar_hmc} that the HMC sample concentrates around the posterior mean of the PDE solution, which is close to the true solution. 
Here, we set $n=100$, $\sigma_e = 0.05$ and generate a set of data using a similar procedure. 
The posterior sample under this setup is shown in Figure \ref{fig:lidar_hmc_sparse}. 
It can be seen that the range of the posterior sample for PDE solution are significantly wider than that in Figure \ref{fig:lidar_hmc}, showcase our method's proficiency in uncertainty quantification with smaller datasets.}

\subsection{Burgers’ Equation - A Nonlinear PDE}\label{subsec:eg_burger_eq}
We consider the viscous Burgers’ equation given in (\ref{eq:burgersequation}), where $s \in [0, 1], t\in[0,0.1]$, with the initial condition $u(0,s) = \exp(-100\times(s-0.5)^2), s\in [0,1]$ and boundary condition $u(t,0)  =0$ and $u(t,1)  =0$, $t\in[0,0.1]$. 
In this PDE, $\theta_2 > 0$ is the viscosity of the fluid, and $\theta_1$ accounts for the scale of the solution. In some research $\theta_1$ is assumed known and equal to $1$. 
In this section, in order to evaluate the performance of parameter estimation of the proposed method and compare with results in \cite{liu2021automated}, we assume that $\theta_1$ is unknown. 
As stated in Section \ref{subsec:PIGPI2}, the Burgers' equation is a nonlinear PDE, and augmentation is employed in this example.
The true value for $\bm{\theta} = (\theta_1,\theta_2)$ is $\bm{\theta}_0 = (1,0.1)$.
We solve the PDE using MATLAB build-in function \textit{pdpde} on a $40000\times 40000$ grid. 
This configuration requires 608 seconds to solve the PDE once. 
Meanwhile, solving this PDE requires a large amount of RAM (Almost 32GB of memory usage is monitored).
For MLE, a global optimization requires at least one thousand evaluations of the PDE solution to guarantee reliable results. 
For BOM, at least $20p = 40$ repeated evaluations are required for each case. 
Each case is repeated 1000 times in this example. 
Thus, methods involving repeated evaluation of numerical solutions (i.e., MLE and BOM) are too time-consuming to apply.\par
In this section, the observation data are generated on a $20\times20$ grid $\bm{\tau} = \{(\frac{2i-1}{400},\frac{2i-1}{40}), i=1,\dots,20, j=1,\dots,20\}$. 
The random errors are i.i.d. normal distributed with zero mean and variance $\sigma_e^2$.
Two cases for $\sigma_e$ are considered, i.e., $\sigma_e = 0.01$ and $\sigma_e = 0.001$.\par
In this example, four methods are evaluated, including PIGP with IBCs (PIGP w IBCs), PIGP without IBCs (PIGP w/o IBCs), and two benchmark methods API and TSM. 
Results are shown in Table \ref{tb:burger_compare_api}.\par
\begin{figure}
    \centering
    \includegraphics[scale =0.5 ]{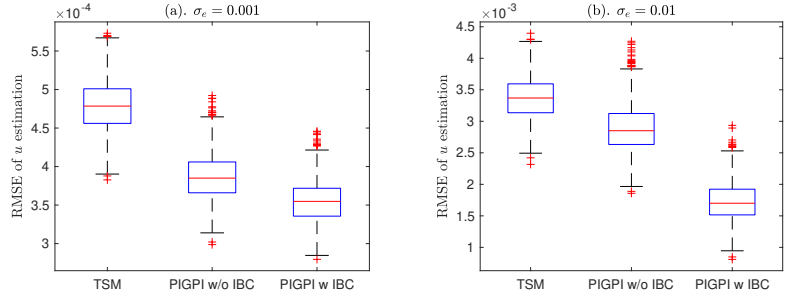}
    \caption{The RMSEs of PDE solution estimation}
    \label{fig:compare_u_est}
\end{figure}

\textbf{Impact of Incorporating IBCs:} In this example, the IBCs are all Dirichlet type. 
To incorporate IBCs, we set $I_1 =(\{0.1,0.2,\dots,0.9\}\times\{0,1\} ) \cup ( \{0\}\times\{0,0.1,\dots,0.9,1\})$.
Except for the case of estimating $\theta_1$ when $\sigma_e = 0.001$, the biases and RMSEs of MAP of parameter $\bm{\theta}$ with IBCs are smaller than the ones without. 
Moreover, Figure \ref{fig:compare_u_est} gives the RMSE of PDE solution estimation on $\bm{I}$. 
It can be seen that the boundary condition has significantly improved the performance of PIGP on the accuracy of PDE solution, where RMSE of $u(\bm{I})$ obtained from PIGP with IBCs are significantly smaller than that from PIGP without IBCs.
Moreover, the PIGP method outperforms the TSM in both cases in terms of estimating the PDE solution.\par
\begin{table}[H]
\caption{The mean of biases, the standard deviation of parameter estimations obtained by PIGP with initial boundary conditions (PIGP w IBC), PIGP method without IBC (PIGP w/o IBC), API method proposed in \cite{liu2021automated} and TSM; coverage rate of $95\%$ credible intervals obtained by two PIGP methods. The proposed methods are emphasized with boldface font. }
\label{tb:burger_compare_api}
\begin{center}
\begin{tabular}{|ll|rr|r|r|}
\hline
\multicolumn{2}{|l|}{\multirow{2}{*}{}}                                                                                         & \multicolumn{2}{l|}{$\sigma_e = 0.001$}                               & \multicolumn{2}{l|}{$\sigma_e = 0.01$}           \\ \cline{3-6} 
\multicolumn{2}{|l|}{}                                                                                                          & \multicolumn{1}{l|}{$\theta_1$}     & \multicolumn{1}{l|}{$\theta_2$} & \multicolumn{1}{l|}{$\theta_1$}        & \multicolumn{1}{l|}{$\theta_2$} \\ \hline
\multicolumn{1}{|l|}{\multirow{4}{*}{\begin{tabular}[c]{@{}l@{}}Bias\\ $\times 10^{-3}$\end{tabular}}} & \textbf{PIGP w IBC}   & \multicolumn{1}{r|}{\textbf{-4.60}} & \textbf{-0.05}                  & \textbf{-15.49}                        & \textbf{-0.27}                  \\ \cline{2-6} 
\multicolumn{1}{|l|}{}                                                                                 & \textbf{PIGP w/o IBC} & \multicolumn{1}{r|}{\textbf{-3.15}} & \textbf{-0.19}                  & \textbf{-23.24}                        & \textbf{-1.60}                  \\ \cline{2-6} 
\multicolumn{1}{|l|}{}                                                                                 & API                    & \multicolumn{1}{r|}{10.76}          & -6.69                           & 108.59                                 & 85.97                           \\ \cline{2-6} 
\multicolumn{1}{|l|}{}                                                                                 & TSM                    & \multicolumn{1}{r|}{-10.44}         & -2.19                           & -50.61                                 & -8.55                           \\ \hline
\multicolumn{1}{|l|}{\multirow{4}{*}{\begin{tabular}[c]{@{}l@{}}SD\\ $\times 10^{-3}$\end{tabular}}}   & \textbf{PIGP w IBC}   & \multicolumn{1}{r|}{\textbf{2.49}}  & \textbf{0.19}                   & \textbf{11.31}                         & \textbf{0.79}                   \\ \cline{2-6} 
\multicolumn{1}{|l|}{}                                                                                 & \textbf{PIGP w/o IBC} & \multicolumn{1}{r|}{\textbf{2.84}}  & \textbf{0.22}                   & \textbf{19.21}                         & \textbf{1.34}                   \\ \cline{2-6} 
\multicolumn{1}{|l|}{}                                                                                 & API                    & \multicolumn{1}{r|}{6.15}           & 0.51                            & 299.73                                 & 258.97                          \\ \cline{2-6} 
\multicolumn{1}{|l|}{}                                                                                 & TSM                    & \multicolumn{1}{r|}{4.01}           & 0.35                            & 27.24                                  & 2.72                            \\ \hline
\multicolumn{1}{|l|}{\multirow{4}{*}{\begin{tabular}[c]{@{}l@{}}RMSE\\ $\times 10^{-3}$\end{tabular}}} & \textbf{PIGP w IBC}   & \multicolumn{1}{r|}{\textbf{5.23}}  & \textbf{0.20}                   & \textbf{19.18}                         & \textbf{0.83}                   \\ \cline{2-6} 
\multicolumn{1}{|l|}{}                                                                                 & \textbf{PIGP w/o IBC} & \multicolumn{1}{r|}{\textbf{4.24}}  & \textbf{0.29}                   & \textbf{30.15}                         & \textbf{2.09}                   \\ \cline{2-6} 
\multicolumn{1}{|l|}{}                                                                                 & API                    & \multicolumn{1}{r|}{12.39}          & 6.71                            & 318.66                                 & 272.74                          \\ \cline{2-6} 
\multicolumn{1}{|l|}{}                                                                                 & TSM                    & \multicolumn{1}{r|}{11.18}          & 2.22                            & 57.46                                  & 8.97                            \\ \hline
\multicolumn{1}{|l|}{\multirow{2}{*}{\begin{tabular}[c]{@{}l@{}}CR\\ $\%$\end{tabular}}}               & \textbf{PIGP w IBC}   & \multicolumn{1}{r|}{\textbf{100}}   & \textbf{100}                    & \textbf{82}                            & \textbf{96.2}                   \\ \cline{2-6} 
\multicolumn{1}{|l|}{}                                                                                 & \textbf{PIGP w/o IBC} & \multicolumn{1}{r|}{\textbf{100}}   & \textbf{100}                    & \textbf{81.4}                          & \textbf{81.5}                   \\ \hline
\end{tabular}
\end{center}
\end{table}

\textbf{Comparing with API proposed in \cite{liu2021automated}:} 
Table \ref{tb:burger_compare_api} shows the biases, standard deviations, and RMSE of MAP obtained by PIGP with and without IBCs, parameter estimation obtained by API method, and parameter estimated using TSM. $\sigma_e$ is specified to be $0.001$ and $0.01$ respectively. 
It can be seen that the proposed method outperforms API in all cases.
From Table \ref{tb:burger_compare_api}, we can see that the RMSE of MAPs obtained by PIGP are orders of magnitude smaller than the RMSE of parameter estimation obtained by API and TSM. 
Another advantage of the PIGP method is the ability to do uncertainty quantification for parameters. 
The normal approximation for posterior density is employed to calculate the $95\%$ credible interval of $\bm{\theta}$ for PIGP. 
The coverage rates of $95\%$ credible intervals are shown in the last two rows of Figure \ref{tb:burger_compare_api}.

\textbf{Implementing HMC and Validating the Accuracy of Normal Approximation:} Now, we validate the accuracy of the normal approximation for the posterior density of parameters by comparing the results with the HMC sample.
In particular, assume the data is observed on a grid $\bm{\tau} = \{(\frac{2i-1}{400},\frac{2i-1}{40}), i=1,\dots,20, j=1,\dots,20\}$. 
The data is polluted by random errors whose distribution is $0$ mean $\sigma_e^2 = 0.01^2$ variance normal distribution.
The initial boundary conditions are included as previously introduced in this section.

\begin{figure}[H]
	\centering \includegraphics[width=15cm]{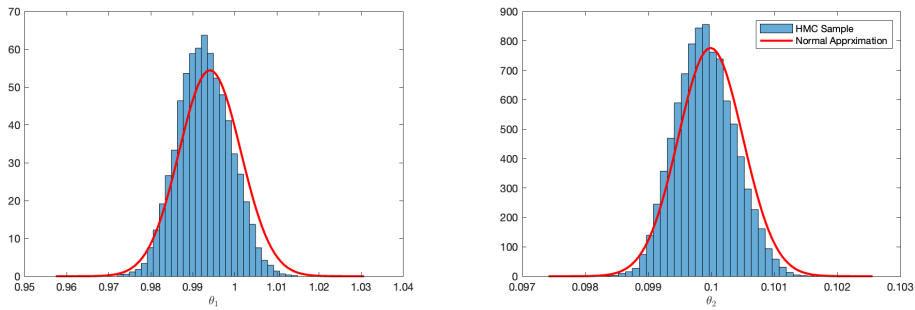}\\
	\centering \caption{Burger's Equation. Posterior sample obtained by HMC, posterior density estimate obtained by normal approximation.}\label{fig:burger_hmc_normal}
\end{figure}

\begin{figure}[H]
	\centering \includegraphics[width=15cm]{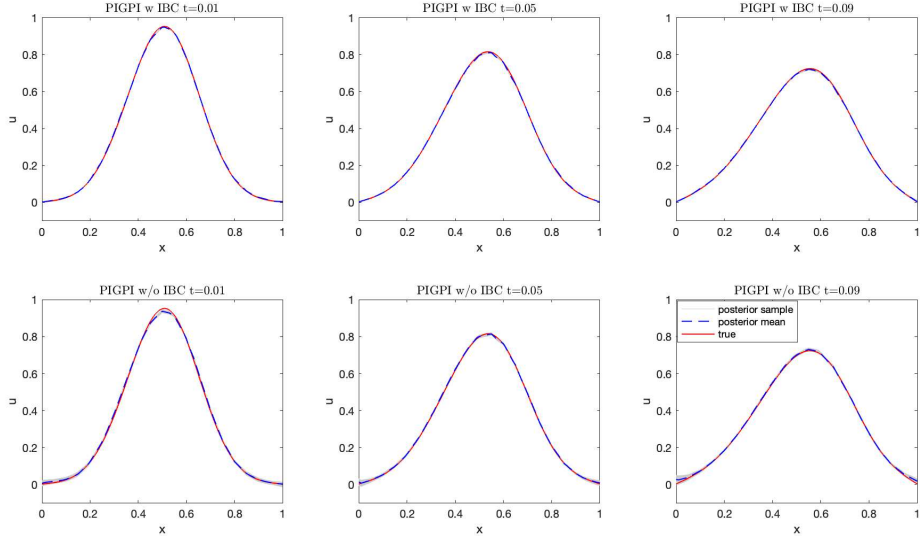}\\
	\centering \caption{Burger's equation. Posterior sample, posterior mean, and the true value of PDE solutions. The posterior sample and posterior mean are produced by HMC using PIGP w/o IBCs and PIGP w IBSs.}\label{fig:hmc_burger}
\end{figure}

In this example, the parameters for HMC are tuned such that the acceptance rate is stable between $60\%$ and $90\%$. 
By comparing the posterior density obtained by HMC and normal approximation, as shown in Figure \ref{fig:burger_hmc_normal}, it can be seen that the normal approximation works well in this example. 
The comparison in this section illustrates the accuracy of the posterior mode and credible interval obtained by normal approximation. 
Figure \ref{fig:hmc_burger} shows the true value, HMC sample, and posterior mean (estimated from HMC sample) of the PDE solution. 
We plot the profiles of the solution at three time values, i.e., $u(0.01,s),u(0.05,s),u(0.09,s), s\in [0,1]$. 
It can be seen that the HMC sample for PIGP w IBCs can cover the true value and the posterior mean is an accurate estimation of the PDE solution for three profiles of the PDE solution. 
As a comparison, the PIGP w/o IBCs produces larger estimation errors than the PIGP w IBCs, especially for $u(t,s)$ when $t = 0.01$, $s$ around $0.5$ and when $t = 0.09$, $s$ around $0$ or $1$. 
The results show that incorporating IBCs can significantly improve the prediction of the PDE solution. 
\subsection{Diffusion-Brusselator Equation - A Nonlinear Multi-Dimensional PDE}\label{subsec:egreaction_diffution}
In this section, we consider a coupled PDE system, i.e., diffusion-Brusselator equation \cite{adomian1995diffusion}.
The PDE system is given by 
\begin{align*}
    \frac{\partial u}{\partial t} & = \theta_1\nabla^2 u + \theta_3 - (\theta_2+1)u+u^2v, &t\times \bm{s}\in [0,1]\times[0,1]^2\\
    \frac{\partial v}{\partial t} & = \theta_1\nabla^2  v+ \theta_2 u -u^2v &t\times \bm{s}\in [0,1]\times[0,1]^2\\
\end{align*}
where $\nabla^2 u = \frac{\partial^2u}{\partial s_1^2}+\frac{\partial^2u}{\partial s_2^2}$ and $\nabla^2 v= \frac{\partial^2v}{\partial s_1^2}+\frac{\partial^2v}{\partial s_2^2}$. 
The boundary conditions are all Neumann boundary conditions describing the adiabatic boundaries, i.e., $\bm{\omega}\cdot \nabla u = 0$, $\bm{\omega}\cdot \nabla v = 0$, $t\in [0,1], \bm{s}\in [0,1]\times\{0,1\}\cup\{0,1\}\times[0,1]$, where $\bm{\omega}$ is the outer normal vector on the specific boundary.  
The initial condition is given by $u(0,\bm{s}) = 2+0.25s_2, v(0,\bm{s}) = 1+0.8 s_1$, $\bm{s}\in [0,1]^2$. The task is to estimate the parameters $\bm{\theta} = (\theta_1, \theta_2, \theta_3)$. 
The augmented PDE is shown in Reaction-Diffusion Systems in Table \ref{tb:aug_example} by setting $F(u,v) =\theta_3 - (\theta_2+1)u+u^2v $ and $G(u,v) = \theta_2 u -u^2v$.
The PDE is solved by MATLAB PDE toolbox, with the mesh generated by setting target maximum edge size \textit{hmax} $= 0.02$, which is slightly smaller than the default setting by MATLAB. 
In average, it takes approximately $17.4$ seconds to solve the PDE once.
In this example, we will examine the performance of the PIGP method on (a) PDE parameter inference for multi-dimensional PDE system in Section \ref{subsubsec:eg_reaction_diffution_full} and (b) PDE parameter inference based on censored data in Section \ref{subsubsec:eg_reaction_diffution_censor}. 

\subsubsection{Parameter Inference for Multi-Dimensional PDE} \label{subsubsec:eg_reaction_diffution_full}
We assume that observation data is given by $y_1(\bm{x}_i) = u(\bm{x}_i) + \varepsilon_{1i}, y_2(\bm{x}_{j+n_1}) = v(\bm{x}_{j+n_1}) + \varepsilon_{2j}, i=1,2,\dots,n_{1},j=1,2,\dots,n_{2}$.
 $\bm{\tau}$ is the union set of $\bm{\tau}_1 = \{\bm{x}_i, i=1,2,\dots,n_{1}\}$ and $\bm{\tau}_2 = \{\bm{x}_{j+n_1},j=1,2,\dots,n_{2}\}$. 
The data are simulated from the PDE solution with true value for $\bm{\theta}$ as  $\bm{\theta}_0 = (0.1,2,1)$.
In this section, $\bm{\tau} = \bm{\tau}_1 = \bm{\tau}_2$ is generated using minimax LHD provided in Python \textit{scipy} package, with specific sample sizes $n = 30, 60, 120, 240$.
The random errors are i.i.d. normal distributed with zero mean and variance $\sigma_e^2 = (0.001^2,0.001^2)$. 
$\bm{I}$ is generated using algorithm \ref{alg:generating_I} with $n_{\bm{I}} = 240$ for each case.\par 
For each $n$, we generate 100 sets of synthetic data and run MLE, BOM, PIGP, and TSM methods on each set of data. 
We compare the bias and standard deviations of parameter estimation obtained from each alternative, the results are shown in table \ref{tb:compare_map_multi_dim}. In this example, MLE is the gold standard that gives the most accurate results. 
However, obtaining MLE involves optimizing the log-likelihood, which requires a large number of evaluations of the PDE numerical solver. 
As shown in Table \ref{tb:compare_map_multi_dim}, the MLE and BOM are much more time-consuming than PIGP\footnote{Adam with 2500 iterations are used for the optimization of the posterior density.} and TSM.
As a surrogate model-based approximation of MLE, BOM performs poorly, especially for the estimation of conductivity $\theta_1$, which is very important in diffusing reaction research. 
For TSM, we can observe that the performance is deteriorating with the decrease in sample size, and TSM performs poorly when the sample size is small for parameter estimation.
In summary, PIGP outperforms BOM and TSM significantly and is the method that is the closest to MLE. 
More comparison results are available in supplementary materials, see Appendix \ref{spsubsec:egreaction_diffution} for details. 
\begin{table}[H]
\caption{The RMSE of parameter inference obtained by MLE, PIGP, TSM, and BOM. The average computational time for parameter inference using MLE, PIGP, TSM, and BOM. The proposed method is emphasized with boldface font. }\label{tb:compare_map_multi_dim}
\begin{center}
\begin{tabular}{|l|l|r|r|r|r|}
\hline
                                                                                          & $n$                & 30            & 60            & 120          & 240          \\ \hline
\multirow{4}{*}{\begin{tabular}[c]{@{}l@{}}$\theta_1$\\ $\times 10^{-3}$\end{tabular}}    & MLE                & 0.3           & 0.2           & 0.2          & 0.1          \\ \cline{2-6} 
                                                                                          & \textbf{PIGP} & \textbf{35.2} & \textbf{10.9} & \textbf{8.3} & \textbf{8.3} \\ \cline{2-6} 
                                                                                          & BOM                & 780.1         & 697.3         & 795.6        & 1044.7       \\ \cline{2-6} 
                                                                                          & TSM                & 169.8         & 58.9          & 21.6         & 13.7         \\ \hline
\multirow{4}{*}{\begin{tabular}[c]{@{}l@{}}$\theta_2$\\ $\times 10^{-3}$\end{tabular}} & MLE                & 0.5           & 0.3           & 0.2          & 0.2          \\ \cline{2-6} 
                                                                                          & \textbf{PIGP} & \textbf{2.7}  & \textbf{1.6}  & \textbf{1.4} & \textbf{1.1} \\ \cline{2-6} 
                                                                                          & BOM                & 87.3          & 88.2          & 86.5         & 92.4         \\ \cline{2-6} 
                                                                                          & TSM                & 1549.1        & 663.7         & 251.4        & 146.3        \\ \hline
\multirow{4}{*}{\begin{tabular}[c]{@{}l@{}}$\theta_3$\\ $\times 10^{-3}$\end{tabular}} & MLE                & 0.5           & 0.4           & 0.3          & 0.2          \\ \cline{2-6} 
                                                                                          & \textbf{PIGP} & \textbf{1.4}  & \textbf{1.1}  & \textbf{1.3} & \textbf{1.2} \\ \cline{2-6} 
                                                                                          & BOM                & 81.0            & 85.0            & 78.9         & 91.2         \\ \cline{2-6} 
                                                                                          & TSM                & 566.4         & 243.2         & 92.2         & 53.7         \\ \hline
\multirow{4}{*}{\begin{tabular}[c]{@{}l@{}}Computational \\ Time\\ (sec)\end{tabular}}    & MLE                & 54328.0       & 55393.2       & 53783.6      & 53258.4      \\ \cline{2-6} 
                                                                                          & \textbf{PIGP}  & \textbf{33.5}            & \textbf{35.3}            & \textbf{35.2}             & \textbf{37.4}            \\ \cline{2-6} 
                                                                                          & BOM                & 1284.7        & 1363.3        & 1383.9       & 1347.3       \\ \cline{2-6} 
                                                                                          & TSM                &12.6                     & 12.9                    & 13.1                     & 15.3                    \\ \hline
\end{tabular}
\end{center}
\end{table}

However, we have discovered that the normal approximation performs poorly in this example, as evidenced by the comparison of the variance of the normal approximation for the posterior with the variance estimated from the HMC samples of the posterior density.
As a result, the normal approximation can significantly underestimate the posterior variances of $\theta_1$ and $\theta_2$, leading to a significant underestimation of the coverage rates of credible intervals. 
Consequently, the normal approximation is not applicable in this case.

\begin{figure}[H]
	\centering \includegraphics[width=15cm]{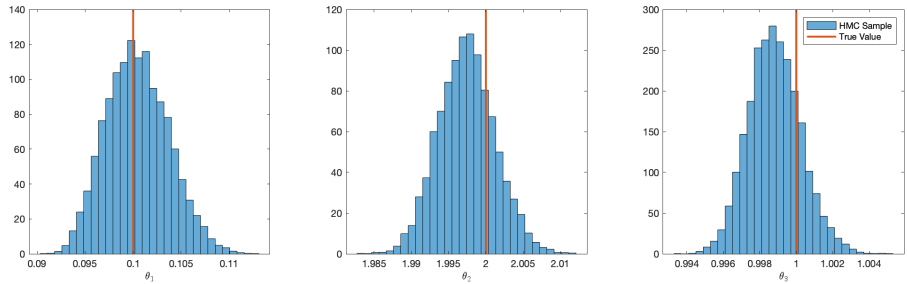}\\
	\centering \caption{The posterior sample of $\theta_1$, $\theta_2$,, and $\theta_3$ obtained by HMC algorithm. Comparing with their true values. }\label{fig:diff_reac_full}
\end{figure}

Given that Table \ref{tb:compare_map_multi_dim} contains 400 simulation cases, we need to run the HMC $400$ times if we expect to calculate the coverage rates  of the confidence interval. 
This process is extremely time-consuming. 
To demonstrate the uncertainty estimation provided by the PIGP method, we selected a specific case, namely when $n=120$.
We generate a set of random data with $n=120$ using the method at the beginning of this subsection.
The HMC algorithm is then applied once to this randomly selected data. 
We utilize the HMC to draw 22,000 random samples, discarding the first $2000$ as burn-in samples. 
The posterior density of the parameter vector $\bm{\theta}$ obtained from the HMC algorithm is reported in Figure \ref{fig:diff_reac_full}.
It can be seen that the posterior samples of $\theta_1$, $\theta_2$, and $\theta_3$ accurately encompass their true values.

\subsubsection{PDE Parameter Inference with Censored Data} \label{subsubsec:eg_reaction_diffution_censor}
In this section, we evaluate the performance of PIGP for estimating $\bm{\theta}$ and PDE solution $u,v$ from partial observations or censored data.
Assume that observation data is given by $y_1(\bm{x}_i) = u(\bm{x}_i) + \varepsilon_{1i}, i=1,2,\dots,n$, i.e., $v$ is fully censored while $u$ is completely missing.
The data are simulated from the PDE solution with true value for $\bm{\theta}$ as  $\bm{\theta}_0 = (0.1,2,1)$. 
In this section, $\bm{\tau}$ is generated using the quasi-Monte Carlo method provided in Python \textit{scipy} package, with specific sample sizes $n = 30, 60, 120, 240$.
The random errors are i.i.d. normal distributed with zero mean and variance $\sigma_e^2 = 0.001^2$. 
$\bm{I}$ is generated using algorithm \ref{alg:generating_I} with $n_{\bm{I}} = 240$ for each case.\par 
For each $n$ case, we generate $100$ sets of synthetic data and run MLE, BOM, and PIGP methods on each set of data.
Note that when one component is censored, TSM is not applicable. 
We compare the bias of parameter estimation obtained from each alternative method, the results are shown in Table \ref{tb:compare_map_censored}.
It can be seen that PIGP outperforms BOM significantly especially for estimation of $\theta_1$. 
As shown in Table \ref{tb:compare_map_multi_dim}, the MLE and BOM are much more time-consuming than PIGP.
We note that in this subsection LBFGS \footnote{LBFGS with 1500 iterations are used for the optimization of the posterior density.} is adopted to solve the MAP optimization with 1500 iterations for PIGP, which can be significantly more time-consuming than the Adam algorithm employed in Section \ref{subsubsec:eg_reaction_diffution_full}. 
However, the computational time of PIGP is still much shorter than MLE and BOM.
Moreover, we can observe that the performances of the proposed method based on censored data are close to that based on full data. 
More comparison results are available in supplementary materials, see Appendix \ref{spsec:simulation_complement} for details. 
\begin{table}[H]
\caption{The RMSE of parameter inference obtained by MLE, PIGP, and BOM.  The average computational time for parameter inference using MLE, PIGP, TSM, and BOM. The proposed method is emphasized with boldface font.}\label{tb:compare_map_censored}    
\begin{center}
\begin{small}
\begin{tabular}{|l|l|r|r|r|r|}
\hline
                                                                                          & $n$                & 30             & 60            & 120           & 240           \\ \hline
\multirow{3}{*}{\begin{tabular}[c]{@{}l@{}}$\theta_1$\\ $\times 10^{-3}$\end{tabular}} & MLE                & 0.2            & 0.2           & 0.2           & 0.1           \\ \cline{2-6} 
                                                                                          & \textbf{PIGP} & \textbf{49.3}  & \textbf{7.0}    & \textbf{5.2}  & \textbf{3.8}  \\ \cline{2-6} 
                                                                                          & BOM                & 1760.3         & 1208.9        & 464.8         & 1385.3        \\ \hline
\multirow{3}{*}{\begin{tabular}[c]{@{}l@{}}$\theta_2$\\ $\times 10^{-3}$\end{tabular}} & MLE                & 0.8            & 0.6           & 0.4           & 0.3           \\ \cline{2-6} 
                                                                                          & \textbf{PIGP} & \textbf{17.9}  & \textbf{3.7}  & \textbf{2.5}  & \textbf{2.4}  \\ \cline{2-6} 
                                                                                          & BOM                & 272.3          & 271.9         & 275.1         & 269.6         \\ \hline
\multirow{3}{*}{\begin{tabular}[c]{@{}l@{}}$\theta_3$\\ $\times 10^{-3}$\end{tabular}} & MLE                & 1.3            & 1.1           & 0.8           & 0.6           \\ \cline{2-6} 
                                                                                          & \textbf{PIGP} & \textbf{321.2} & \textbf{64.9} & \textbf{36.5} & \textbf{31.9} \\ \cline{2-6} 
                                                                                          & BOM                & 499.2          & 487.7         & 504.2         & 488.5         \\ \hline
\multirow{3}{*}{\begin{tabular}[c]{@{}l@{}}Computational \\ Time\\ (sec)\end{tabular}}    & MLE                & 43505.0        & 41217.6       & 45458.9       & 49449.9       \\ \cline{2-6} 
                                                                                          & PIGP            & \textbf{469.1}                    &\textbf{465.4}                    & \textbf{470.2}                    & \textbf{471.8}                    \\ \cline{2-6} 
                                                                                          & BOM                & 1354.0         & 1358.9        & 1369.1        & 1358.4        \\ \hline
\end{tabular}
\end{small}    
\end{center}
\end{table}

Moreover, the PIGP method is able to estimate the missing components. 
Other alternatives such as MLE and BOM all rely on the numerical PDE solution to give the estimation of the PDE solution. 
As a comparison, we estimate $u$ and $v$ simutaneously on a grid $\{0.1,0.3,\dots,0.9\}^3$ using PIGP, MLE and BOM. 
The results are shown in Table \ref{tb:est_u_v_censored}. 
It can be seen that MLE is the best method ignoring the computational time. 
When $n$ is very small, say $n=30$, the PIGP cannot estimate censored components well. 
However, for larger $n$, the PIGP can estimate both $u$ and $v$ more accurately than the BOM, with a significant time saving as shown in Table \ref{tb:compare_map_censored}. 
\begin{table}[H]
\caption{The mean of RMSE of PDE solution estimation, $v$ is censored. The PIGP is able to reconstruct the fully censored components. The proposed method is emphasized with boldface font.}\label{tb:est_u_v_censored}
\begin{center}
\begin{tabular}{|ll|r|r|r|r|}
\hline
\multicolumn{2}{|l|}{$n$}                                                                                                          & 30              & 60             & 120            & 240           \\ \hline
\multicolumn{1}{|l|}{\multirow{3}{*}{\begin{tabular}[c]{@{}l@{}}u\\ $\times 10^{-3}$\end{tabular}}}            & MLE            & 0.23            & 0.18           & 0.13           & 0.10          \\ \cline{2-6} 
\multicolumn{1}{|l|}{}                                                                                            & \textbf{PIGP} & \textbf{10.35}  & \textbf{2.81}  & \textbf{1.54}  & \textbf{1.25} \\ \cline{2-6} 
\multicolumn{1}{|l|}{}                                                                                            & BOM            & 58.55           & 61.56          & 58.69          & 60.86         \\ \hline
\multicolumn{1}{|l|}{\multirow{3}{*}{\begin{tabular}[c]{@{}l@{}}v(censored) \\ $\times 10^{-3}$\end{tabular}}} & MLE            & 0.34            & 0.31           & 0.21           & 0.16          \\ \cline{2-6} 
\multicolumn{1}{|l|}{}                                                                                            & \textbf{PIGP} & \textbf{159.81} & \textbf{26.95} & \textbf{11.78} & \textbf{8.57} \\ \cline{2-6} 
\multicolumn{1}{|l|}{}                                                                                            & BOM            & 92.12           & 101.70         & 102.14         & 97.15         \\ \hline
\end{tabular}
\end{center}
\end{table}

\begin{figure}[H]
	\centering \includegraphics[width=15cm]{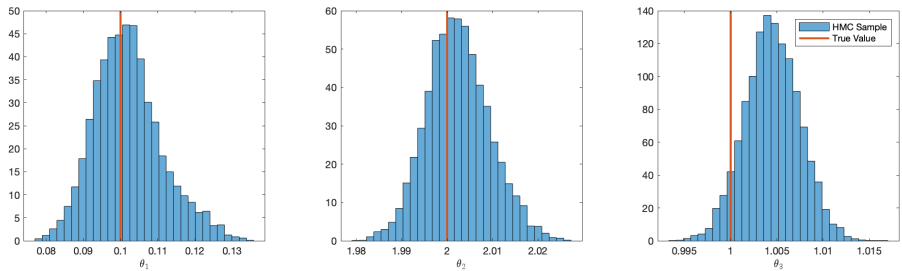}\\
	\centering \caption{The posterior sample of $\theta_1$, $\theta_2$,, and $\theta_3$ obtained by HMC algorithm. Compared with their true values. }\label{fig:diff_reac}
\end{figure}

The same phenomenon is observed in the case of censored data: we found that the normal approximation significantly underestimates the posterior variances for $\bm{\theta}$.
As an analogy to Section \ref{subsubsec:eg_reaction_diffution_censor}, a random set of data (with $v$ censored) is generated with $n=240$ using the method mentioned at the beginning of this subsection.  
The HMC algorithm is applied to this randomly selected data. 
We use the HMC to draw $22,000$ random samples, discarding the first $2,000$ as burn-in samples. 
The resulting posterior density of the parameter vector $\bm{\theta}$ is reported in Figure \ref{fig:diff_reac}. 
The results indicate that the posterior samples of $\theta_1$, $\theta_2$, and $\theta_3$ accurately capture their true values.

An interesting finding arises when comparing Figure \ref{fig:diff_reac_full} and Figure \ref{fig:diff_reac}. 
In this subsection, we compare the results of the HMC sample with those presented in Section \ref{subsubsec:eg_reaction_diffution_full}. 
In both cases, the same amount of data is used, i.e., in Section \ref{subsubsec:eg_reaction_diffution_full}, both components have $120$ observed data points, while in this section, only $u$ is observed at 240 points, and $v$ is completely censored. 
By comparing Figure \ref{fig:diff_reac_full} and \ref{fig:diff_reac}, it is evident that the posterior variances of $\theta_1$, $\theta_2$, and $\theta_3$ obtained from the full data (Section \ref{subsubsec:eg_reaction_diffution_full}) are significantly smaller than the posterior variances obtained from the censored data (as shown in this subsection), given that both posterior samples can cover the true value of $\bm{\theta}$ well. 
Thus, even using the same amount of data, the full data provide more information than the censored data, resulting in smaller posterior variances for the unknown parameters. 

\section{Conclusion and Discussions}\label{sec:conclusion}
In this paper, we propose a novel methodology, called the PIGP method, for PDE parameter inference. 
We assign a GP prior on the PDE solution.
By constraining the GP on a manifold that compiles with the PDE and its boundary and initial conditions, the posterior density of unknown parameters and PDE solution is constructed with the observation data and manifold constraint. 
To handle nonlinear and parameter-dependent PDE operators, we propose an augmentation method to construct an equivalent PDE system with zeroth order nonlinearity.
Two enhancements are provided to ensure computational efficiency.
First, we propose a simple algorithm to select discretization point set $\bm{I}$ to achieve good space-filling properties. 
Second, we propose to use KL expansion to GP that can provide a set of independent parameters. 
Moreover, the truncated KL expansion can reduce the dimension of parameter space significantly.
HMC and normal approximation are adopted for the simultaneous inference and uncertainty quantification of parameters and PDE solutions. \blue{Despite its advantages, PIGP could have the limitation that its performance depends on the quality of the experimental design related to the discretization point set $\bm{I}$. Furthermore, in cases where extrapolation is required, 
$\bm{I}$ must cover the extrapolation region. These considerations represent common challenges in the realm of surrogate modeling \cite{sacks1989design, queipo2005surrogate}. }

We employ four examples to illustrate the performance of the PIGP method on a broad range of applications. 
First, we consider a simple toy example to investigate the impact of the choice of discretization set. 
We discover that when observation data is sparse, choosing a large discretization set helps improve the accuracy of parameter estimation. 
In this example, we also demonstrate the performance of the PIGP method when IBCs are unavailable.
Second, the LIDAR example is employed to illustrate the performance of the PIGP method by comparing it with methods proposed in \cite{xun2013parameter} and TSM.
In the third example, we employ a nonlinear example (Burger equation) to compare the performance of the PIGP method with the API method proposed in \cite{liu2021automated} and TSM.
In the last example, we employ a multi-dimensional nonlinear PDE to demonstrate the ability of PIGP to handle multi-output systems and missing data. 
In summary, we illustrate the PIGP method using two criteria: accuracy and efficiency.

The PIGP method is not perfect.
There are still some limitations for PIGP that need further investigation. 
One potential example is the heat moisture coupled PDE, in which two components represent the moisture content and temperature respectively. 
The values and ranges of temperature can be $100$ times larger than the moisture content. 
Thus, when moisture data is observable, the PIGP can be applied to infer parameters from temperature data. 
However, the hyperparameters of GP prior for the moisture should be very different from those for temperature. 
In this case, physical knowledge can be used for initializing the GP prior for the missing components. 
Second, if PDE is very complicated, eg. the PDE operator consists of many derivative terms ($l$ is large, see Section \ref{subsec:PIGPI2} for definition) that involve high order derivatives of the solution, one needs to choose a GP prior with a large value for the degree of smoothness in the Mat\'{e}rn kernel. 
On the other hand, when $l$ is large and augmentation is required, the augmented PDE will have many components. 
Both situations can lead to a nearly singular covariance matrix for the GP and its derivatives ($C$ and $\mathcal{K}$ in \eqref{eq:general_posterior}).
We will continue investigating how to improve the numerical stability of PIGP when PDE is very complicated.
The third limitation is that the normal approximation is not applicable in some situations. 
For example, in Section \ref{subsec:egreaction_diffution}, we found that the normal approximation will suggest a significantly smaller posterior variance compared with the variance estimated from the HMC sample. 
Thus, a more robust and efficient approximation algorithm for posterior uncertainty estimation is an important future work.


There are potential directions worth further investigation.
First, when physical data is insufficient, say $n=10$ or even unavailable, one might need to do follow-up physical experiments. 
The design of physical experiments in a sequential manner is an important task \cite{wu2021experiments}. 
Second, we focus on the inference of scalar/vector-valued parameters in this paper.
When the parameter of interest is a function instead of a constant\cite{li2021gaussian}, estimating the functional parameter can be an important extension of the current method. 

\bibliographystyle{siam}
\bibliography{ref_for_magi.bib} 

\newpage

\appendix
\section{Gaussian Processes}\label{spsec:gpr}
All details about the Gaussian process(GP) are in this section.
In Section \ref{spsubsec:kernel_choice}, we discuss the choice of kernel functions for GP priors of PDE.
Section \ref{spsubsec:matern} provides the detailed derivatives for the covariance matrices that appear in the paper. 
Section \ref{spsubsec:train_GP} discusses the initialization of GP models. 
Section \ref{spsubsec:predict_GP} introduces the prediction of $u(\bm{x})$ where $\bm{x}\notin \bm{I}$.
\subsection{Choices of Kernel Functions}\label{spsubsec:kernel_choice}
For GP with multi-dimensional input variable $\bm{x}\in \mathbb{R}^p$, where $p> 1$, two types of Mat\'{e}rn kernel are commonly employed as the correlation functions.  
The first type is the isotropic Mat\'{e}rn kernel, which has the form $\mathcal{K}(\bm{x},\bm{x}')=\psi_1 \frac{2^{1-\nu}}{\Gamma(\nu)}(\sqrt{2\nu}\frac{\|\bm{x}-\bm{x}'\|}{\psi_2})^{\nu}B_{\nu}(\sqrt{2\nu}\frac{\|\bm{x}-\bm{x}'\|}{\psi_2})$. 
The second type is the product Mat\'{e}rn kernel, which has the form $\mathcal{K}(\bm{x},\bm{x}')=\psi_1\prod_{i=1}^p \frac{2^{1-\nu}}{\Gamma(\nu)}(\sqrt{2\nu}\frac{|x_i-x_i'|}{\psi_{2i}})^{\nu}B_{\nu}(\sqrt{2\nu}\frac{|x_i-x_i'|}{\psi_{2i}})$. 
In both types of kernels, $\Gamma$ is the Gamma function, $B_{\nu}$ is the modified Bessel function of the second kind, and $\nu$ is the degree of smoothness parameter.\par
In this paper, we choose the product Mat\'{e}rn kernel.
There are two reasons for our preference.
First, the product Mat\'{e}rn kernel has the advantage over isotropic Mat\'{e}rn kernel in that the smoothness (i.e., the length-scale $\psi_{2i}$) of $u$ to different $x_i$ can be different.
Thus, the product Mat\'{e}rn kernel is able of tuning the variability of the GP in each direction.
Second, the second-order derivative of isotropic Mat\'{e}rn kernel is not well-defined at the origin.  
A well-defined second-order derivative of the kernel function is important to ensure the Gaussianity of the derivative of the original GP. 
Thus, the GP with isotropic Mat\'{e}rn kernel is not suitable for PIGP. 
The result is summarized as proposition \ref{prop:iso_matern_prop}. 
\begin{proposition}\label{prop:iso_matern_prop}
Suppose the kernel function is isotropic Mat\'{e}rn kernel with $\nu>2$. 
Then the second order derivative of $\mathcal{K}(\bm{x},\bm{x}')$ with respect to $x_i$ and $x'_i$ is not continuous at $\bm{x}=\bm{x}'$ for any $i=1,\dots,p$.
\end{proposition}
\begin{proof}
Consider 
\begin{align*}
  	    \frac{\partial^2 \mathcal{K}(\bm{x},\bm{x}')}{\partial x_1 \partial x_1'} &= \frac{\partial^2 \mathcal{K}_{\nu}(r)}{\partial r^2}\left(\frac{\partial r}{\partial x_1}\frac{\partial r}{\partial x_1'}\right)+\frac{\partial \mathcal{K}_{\nu}(r)}{\partial r}\left(\frac{\partial^2 r}{\partial x_1\partial x_1'}\right)\\
  	   & = -\frac{\partial^2 \mathcal{K}_{\nu}(r)}{\partial r^2}\left(\frac{(x_1-x_1')^2}{r^2}\right)+\frac{\partial \mathcal{K}_{\nu}(r)}{\partial r}\left(\frac{3(x_1-x_1')^2}{r^3}-\frac{1}{r}\right)
\end{align*}
Since $\frac{\partial \mathcal{K}_{\nu}(r)}{\partial r}\frac{1}{r}$ and $\frac{\partial^2 \mathcal{K}_{\nu}(r)}{\partial r^2}$ are well defined and nonzero(see Subsection \ref{spsubsec:matern}) for $r\geq0$ and limit when $r\to 0$ exists.
We consider the value $\frac{(x_1-x_1')^2}{r^2}$. 
Assume $\bm{x}= (x_1,x_2), \bm{x}' = (x_1',x_2')$ are 2-d variables. 
Fixing $\bm{x}$, consider the limit of $\frac{(x_1-x_1')^2}{r^2}$ when $\bm{x}'\to \bm{x}$ from different direction.
We first assume that $x_2=x_2' = 0$, i.e., let $\bm{x}'\to \bm{x}$ from $x_1$ axis.
The limit will be $\frac{(x_1-x_1')^2}{r^2}\to 1$. 
Then we assume that $x_2'-x_2 = x_1'-x_1$. 
The limit will be $\frac{(x_1-x_1')^2}{r^2}\to \frac{1}{2}$. 
From these values, we see that there is no proper definition for $\frac{(x_1-x_1')^2}{r^2}$ when $r=0$ (cannot be defined by limit). 
\end{proof}
In this paper, we recommend $\nu = 2a+0.1$ for all kernels, where $a = \max_{1\leq i\leq p}{\|\bm{\alpha}_i\|_{\infty}}$ to ensure that the derivatives of GP provide continuous differential properties up to $2\alpha$ order for each coordinate $x_1$. 
This property is important since it ensures that the definition $\mathcal{LK}, \mathcal{KL}$, and $\mathcal{LKL}$ exists. 
\subsection{Derivatives for Product Mat\'{e}rn Kernel}\label{spsubsec:matern}
For Product Mat\'{e}rn kernel, we have $\mathcal{K}(\bm{x},\bm{x}')=\prod_{i=1}^p \mathcal{K}_i(x_i,x_i')=\prod_{i=1}^p \mathcal{K}_{\nu}(r_i)$, where $r_i=\sqrt{2v}\frac{|x_i-x_i'|}{\psi_{2i}}$.	 
To calculate the derivatives of $\mathcal{K}(\bm{x},\bm{x}')$, i.e., $\mathcal{L}_{\bm{x}}^{\bm{\theta}}\mathcal{K}(\bm{x},\bm{x}')$, $\mathcal{L}_{\bm{x}'}^{\bm{\theta}}\mathcal{K}(\bm{x},\bm{x}')$ and $\mathcal{L}_{\bm{x}'}^{\bm{\theta}}\mathcal{L}_{\bm{x}}^{\bm{\theta}}\mathcal{K}(\bm{x},\bm{x}')$, it suffices to calculate the quantity:
\begin{equation}\label{eq:derivative_k0}
	\frac{\partial^{|\bm{\alpha}|} \mathcal{K}(\bm{x},\bm{x}')}{\partial x_1^{\alpha_1}\cdots x_p^{\alpha_p}}, \frac{\partial^{|\bm{\alpha}|} \mathcal{K}(\bm{x},\bm{x}')}{\partial {x'}_1^{\alpha_1}\cdots {x'}_p^{\alpha_p}}, \frac{\partial^{|\bm{\alpha}|+|\bm{\alpha}'|} \mathcal{K}(\bm{x},\bm{x}')}{\partial x_1^{\alpha_1}\cdots x_p^{\alpha_p}{x'}_1^{\alpha_1'}\cdots {x'}_p^{\alpha_p'}},
\end{equation}
for all $\bm{\alpha}\in A$. 
Due to the definition of $r_i$, $\frac{\partial r_i}{\partial x_i}=\frac{\partial r_i}{\partial d_i}\frac{\partial d_i}{\partial x_i}=\frac{\sqrt{2v}}{\psi_2} \frac{x_i-x_i'}{|x_i-x_i'|}$, we have $\frac{\partial r_j}{\partial x_i}=0$ for $j\neq i$, and $\frac{\partial^2 r_i}{\partial x_i^2}=0$. 
Then,
\begin{equation}\label{eq:derivative_k1}
	\frac{\partial^{|\bm{\alpha}|} \mathcal{K}(\bm{x},\bm{x}')}{\partial x_1^{\alpha_1}\cdots x_p^{\alpha_p}}=\prod_{i=1}^p \frac{\partial^{\alpha_i} \mathcal{K}_{\nu}(r_i)}{\partial r_i^{\alpha_i}}\left(\frac{\partial r_i}{\partial x_i}\right)^{\alpha_i},
\end{equation}
and
\begin{equation}\label{eq:derivative_k2}
	\frac{\partial^{|\bm{\alpha}|+|\bm{\alpha}'|} \mathcal{K}(\bm{x},\bm{x}')}{\partial x_1^{\alpha_1}\cdots x_p^{\alpha_p}{x'}_1^{\alpha_1'}\cdots {x'}_p^{\alpha_p'}}=\prod_{i=1}^p \frac{\partial^{\alpha_i+\alpha_i'} \mathcal{K}_{\nu}(r_i)}{\partial r_i^{\alpha_i+\alpha_i'}}\left(\frac{\partial r_i}{\partial x_i}\right)^{\alpha_i}\left(\frac{\partial r_i}{\partial x_i'}\right)^{\alpha_i'}. 
\end{equation}
The key step is to calculate $\frac{\partial^{\alpha_i} \mathcal{K}_{\nu}(r_i)}{\partial r_i^{\alpha_i}}$ and $\frac{\partial^{\alpha_i+\alpha_i'} \mathcal{K}_{\nu}(r_i)}{\partial r_i^{\alpha_i+\alpha_i'}}$. 
It reduces to the calculation of the derivatives of 1d Mat\'{e}rn kernel. 
For notaional convenience, let $\mathcal{K}_{\nu}(r)=\frac{2^{1-\nu}}{\Gamma(\nu)}r^{\nu}B_{\nu}(r)$ and $r\geq 0$. 
We herein calculate $\mathcal{K}_{\nu}(r)$ up to $2a$ order. 
Define $A_{\nu}(r) = r^{\nu}B_{\nu}(r)$
\begin{equation*}
	\mathcal{K}_{\nu}(r)=\frac{2^{1-\nu}}{\Gamma(\nu)}r^{\nu}B_{\nu}(r)=\frac{2^{1-\nu}}{\Gamma(\nu)}A_{\nu}(r),
\end{equation*}
then,
\begin{align*}
	\frac{\partial^i \mathcal{K}_{\nu}(r)}{\partial r^i}&=	\frac{2^{1-\nu}}{\Gamma(\nu)}\frac{\partial^i A_{\nu}(r)}{\partial r^i}. 
\end{align*}
Since
\begin{align*}
	\frac{\partial A_{\nu}(r)}{\partial r}&=v r^{v-1}B_{\nu}(r)+r^{\nu} \frac{\partial B_{\nu}(r)}{\partial r}\\
	&=v r^{v-1}B_{\nu}(r)-r^{\nu} \frac{B_{v-1}(r)+B_{v+1}(r)}{2}\\
	&=\frac{v A_{\nu}(r)}{r}-\frac{rA_{v-1}(r)}{2}-\frac{A_{v+1}(r)}{2r}\\
	&=\frac{v A_{\nu}(r)}{r}-\frac{rA_{v-1}(r)}{2}-\frac{2v A_{\nu}(r)+r^2 A_{v-1}(r)}{2r}\\
	&=-r A_{v-1}(r),
\end{align*}
given that $A_{v+1}(r)=r^{v+1}B_{v+1}(r)=r^{v+1}(\frac{2vB_{\nu}(r)}{r}+B_{v-1}(r))=2v A_{\nu}(r)+r^2 A_{v-1}(r)$. 
Repeating this process, we have 
\begin{align*}
	\frac{\partial^2 A_{\nu}(r)}{\partial r^2} &=-\frac{\partial rA_{v-1}(r)}{\partial r}\\
	&=-A_{v-1}(r)+r^2 A_{v-2}(r);
\end{align*}
\begin{align*}
	\frac{\partial^3 A_{\nu}(r)}{\partial r^3} &=\frac{\partial (-A_{v-1}(r)+r^2 A_{v-2}(r))}{\partial r}\\
	&=3rA_{v-2}(r)-r^3 A_{v-3}(r);
\end{align*}
\begin{align*}
	\frac{\partial^4 A_{\nu}(r)}{\partial r^4} &=\frac{\partial (3rA_{v-2}(r)-r^3 A_{v-3}(r))}{\partial r}\\
	&=3A_{v-2}(r)-6r^2 A_{v-3}(r)+r^4A_{v-4}(r).
\end{align*}
In general, we can derive the differentiation up to $2a$ order by 
\begin{align}\label{eq:derivative_k3}
    \frac{\partial^k A_{\nu}(r)}{\partial r^k} &= \frac{\partial^{k-1} -rA_{v-1}(r)}{\partial r^{k-1}}\nonumber\\
    &= -r\frac{\partial^{k-1} A_{v-1}(r)}{\partial r^{k-1}}-\frac{\partial^{k-2} A_{v-1}(r)}{\partial r^{k-2}},
\end{align}
which can be easily derived from induction.
It is worth noting that the derivation of $\frac{\partial^k \mathcal{K}_{\nu}(r)}{\partial r^k}$ at $r=0$ is well defined because $A_{\nu}(r)$ is well defined and finite on $[0,\infty)$.
Then \eqref{eq:derivative_k0} is immediately obtained by substituting \eqref{eq:derivative_k3} into \eqref{eq:derivative_k1} and \eqref{eq:derivative_k2}

\subsection{Initialization of GP Models for \texorpdfstring{$u(\bm{x})$}{u(x)} and  \texorpdfstring{$u_i(\bm{x})$}{ui(x)}}\label{spsubsec:train_GP}
This section introduces how to determine the hyper-parameters for the prior GP models for $u(\bm{x})$ and $u_i(\bm{x})$. \par
First, the GP model $U(\bm{x})\sim \text{GP}(\mu,\sigma^2 \mathcal{K}(\bm{x},\bm{x}'))$ for $u(\bm{x})$ is trained using observation data, i.e., $y_i = u(\bm{x}_i)+\varepsilon_i$ if they are available. 
Note that for the situation where variance for $\epsilon$ is unknown, we estimate the hyper-parameters of GP and random error variance simultaneously.\par
Second, when PDE augmentation is used, the data for $u_i(\bm{x}), i>1$  are unavailable.
To train the GP model for $U_i(\bm{x}), i>1$, we first generate a set of synthetic data. 
Then the GP model is trained using the synthetic data. 
In particular, suppose $U_1(\bm{I})$ is trained by $y_i(\bm{x}_i )= u_1(\bm{x}_i)+\varepsilon_i$, then the synthetic-data for $U_k = \nabla^{\bm{\alpha}_k} U_1$ can be generated from the derivative of the GP $U_1$, i.e., $\hat{U}_k(\bm{I}) = \nabla^{\bm{\alpha}_k}\mu + \nabla^{\bm{\alpha}_k}_1\mathcal{K}(\bm{I},\bm{I})\mathcal{K}(\bm{I},\bm{I})^{-1} (U_1(\bm{I})-\mu) $, where $U_1(\bm{I})$ is the estimation of $u_1(\bm{I})$ by GP model, the subscript $1$ for $\nabla$ means that we take the derivative of $\mathcal{K}$ with respect to it's first variable. 
Moreover, we assume that the GP priors $U_i,i=1,\dots,l$ are independent.\par
Third, for completely missing components in multi-dimensional PDE problems, we set an initial value for GP hyper-parameter of missing components as equal to the hyper-parameters of GP for observed components.
Then, the hyper-parameters can be set as unknown and optimized in the posterior density optimization loop. 
In Section \ref{subsec:egreaction_diffution}, since the two components have similar variations and magnitudes, we fix the GP hyper-parameter of missing components as equal to the hyper-parameters of GP for the observed component.
However, if variations and magnitudes of two components (missing component and available component) of the PDE are significantly different, the hyper-parameters for the GP model for missing components can be set as unknown hyper-parameters and optimized with parameters $\bm{\theta}$ and $u(\bm{I})$ in the posterior density optimization loop. 
We mention this issue in Section \ref{sec:conclusion}.

\subsection{Prediction of PDE Solution \texorpdfstring{$u(\bm{x})$}{u(x)} for \texorpdfstring{$\bm{x}\notin \bm{I}$}{u(x) notin I} }\label{spsubsec:predict_GP}
Assume $\bm{I}$ is known and fixed.
In this section, we introduce the prediction/estimation of $u(\bm{x})$ when $\bm{x}\notin\bm{I}$.
Estimation of $u(\bm{x})$ is obtained by the conditional mean of GP, i.e., $U(\bm{x})|U(\bm{I}) = \hat{u}(\bm{I})$. 
In particular, $U(\bm{x})|U(\bm{I}) = \hat{u}(\bm{I})\sim N(\mu +\mathcal{K}(\bm{x},\bm{I})\mathcal{K}(\bm{I},\bm{I})^{-1}(\hat{u}(\bm{I})-\mu),\sigma^2(1-\mathcal{K}(\bm{x},\bm{I})\mathcal{K}(\bm{I},\bm{I})^{-1}\mathcal{K}(\bm{I},\bm{x})))$. 
Thus, $m(\bm{x}) = \mu +\mathcal{K}(\bm{x},\bm{I})\mathcal{K}(\bm{I},\bm{I})^{-1}(\hat{u}(\bm{I})-\mu)$ is adopted as the estimation for PDE solution at $\bm{x}$. 
 
\section{Supplementary Discussions on Algorithms}\label{spsec:other_discussion}
This section provides necessary discussions on the implementation of the algorithm proposed in Section \ref{sec:algorithm}, including KL expansion (discussed in Subsection \ref{spsubsec:kl_expansion_remark}), HMC (discussed in Subsection \ref{spsubsec:uq_posterior_remark}) and algorithms for optimization of posterior density (discussed in Subsection \ref{spsubsec:map_opt_alg}). 
\subsection{Remarks on KL-expansion}\label{spsubsec:kl_expansion_remark}
By using KL approximation, we can rewrite the posterior density (\ref{eq:general_posterior}) to 
\begin{align}\label{eq:boundary_modified_posterior_KL}
	&p_{\bm{\Theta}, \bm{Z},\sigma_e^2| W_{\bm{I}}, Y(\bm{\tau})=y(\bm{\tau}), U(\bm{I})=u(\bm{I})}\left(\bm{\theta}, \bm{z},\sigma_e^2|W_{\bm{I}}=0, Y(\bm{\tau})=y(\bm{\tau}), U(\bm{I})=u(\bm{I})\right)\notag \\
	\propto& \exp \Big\{-\frac{1}{2}\big[\left\| \bm{z}\right\|_2^2+(n+2)\log(\sigma_e^2)+\left\| u(\bm{\tau})-y(\bm{\tau})\right\|_{\sigma_e^{-2}}\notag\\ &+\log|K|+\left\|f(\bm{I},u(\bm{I}),\bm{\theta})-\mathcal{L}_{\bm{x}}^{\bm{\theta}}\mu(\bm{I})-m\{u(\bm{I})-\mu(\bm{I})\} \right\|_{K^{-1}}\big]  \Big\},
\end{align}
where $u(\bm{I})=\sum_{i=1}^M z_i \sqrt{\lambda_i}\bm{\varphi}_i$ and $u(\bm{\tau})$ is extracted from $u(\bm{I})$. 
To further speeding up the computation of evaluating (\ref{eq:boundary_modified_posterior_KL}), we can use
\begin{align*}
   & \left\|f(\bm{I},u(\bm{I}),\bm{\theta})-\mathcal{L}_{\bm{x}}^{\bm{\theta}}\mu(\bm{I})-m\{u(\bm{I})-\mu(\bm{I})\}\right\|_{K^{-1}}\\
   =& \left\|f(\bm{I},u(\bm{I}),\bm{\theta})-\mathcal{L}_{\bm{x}}^{\bm{\theta}}\mu(\bm{I})-m\{\bm{\Phi}_M\bm{z}\} \right\|_{K^{-1}} \\
   = &  \left\|\bm{z}- \bm{m}_z \right\|_{\bm{\Phi}_M^Tm^TK^{-1}m\bm{\Phi}_M},
\end{align*}
where $\bm{m}_z =(\bm{\Phi}_M^T m_u^Tm_u \bm{\Phi}_M)^{-1}(\bm{\Phi}_M^T m_u^T)(f(\bm{I},u(\bm{I}),\bm{\theta})-\mathcal{L}_{\bm{x}}^{\bm{\theta}}\mu(\bm{I}))$.
By discussion in previous sections, $K_c$ can be pre-calculated and stored, now we can pre-calculate $\bm{\Phi}_M^Tm^TK^{-1}m\bm{\Phi}_M$. 
Thus the computational cost of posterior density has reduced to $O(n_{\bm{I}}M)$, which is dominated by computation of $\bm{m}_z$.
\subsection{Hamiltonian Monte Carlo}\label{spsubsec:uq_posterior_remark}
We apply the Hamiltonian Monte Carlo(HMC) algorithm to draw random samples from the posterior distribution. 
A more thorough introduction to HMC is available in [34].
For simplicity, let $\bm{q}$ denote the collection of $U(\bm{I})$, $\bm{\theta}$ and $\sigma_e^2$ (if assumed unknown). 
Suppose that the posterior density is $p(\bm{q}) = \frac{1}{C}\exp\{-L(\bm{q})\}$, where $K$ is a normalized constant, $L(\bm{q})$ is the minus log posterior density.
Then we can introduce the momentum variables, $\bm{p}$, which has the same dimension as $\bm{q}$. 
Define the "Kinetic energy" to be $K(\bm{p}) = \bm{p}^T\bm{p}/2$. 
Then the Hamiltanian is defined to be $H(\bm{q},\bm{p})) = L(\bm{q})+K(\bm{p})$. 
Consider the joint density of $\bm{q}$ and $\bm{p}$, which is proportional to $\exp(-H(\bm{q}, \bm{p}))$. 
Under this construction, $\bm{q}$ and $\bm{p}$ are independent, where the marginal probability density of $\bm{q}$ is the target posterior density and the marginal probability density of $\bm{p}$ is Gaussian. 
We will then sample from this augmented distribution for $(\bm{q}, \bm{p})$. 
We repeat the following three steps, which together compose one HMC iteration: 
(1). Sample $\bm{p}$ from the normal distribution $N(0, I)$ since $K(\bm{p}) = \bm{p}^T\bm{p}/2$ corresponding to standard normal distribution; 
(2). construct a proposal $(\bm{q}_*, \bm{p}_*)$ for $(\bm{q}, \bm{p})$ by simulating the Hamiltonian dynamics using the leapfrog method (detailed in the next paragraph), and 
(3). accept or reject $(\bm{q}_*, \bm{p}_*)$ as the next state of $(\bm{q}, \bm{p})$ according to the Metropolis acceptance probability, i.e.,$\min\{1, \exp(-H(\bm{q}_*, \bm{p}_*) +H(\bm{q}, \bm{p}))\}$.
After repeating the HMC iteration for the desired number of iterations, the sampled $q$ are taken to be the samples from posterior density. \\
\textbf{Leapfrog Method for Hamiltonian Dynamics}. 
The leapfrog method is used to approximate the Hamiltonian dynamics.
One step of the leapfrog method with step size $\varepsilon$ from an initial point $(\bm{q}_0, \bm{p}_0)$ consists of three parts. First, we make a half step for the momentum, $\bm{p} = \bm{p}_0- (\varepsilon/2)\nabla U(\bm{q})|\bm{q}=\bm{q}_0$.
Second, we make a full step for the position, $\bm{q}_* = \bm{q}_0 + \varepsilon\bm{p}$. 
Third, we make another half step for the momentum using the gradient evaluated at the new position, $\bm{p}_* = \bm{p}-(\varepsilon/2)\nabla U(\bm{q})|\bm{q}=\bm{q}_*$.
The step size $\varepsilon$ and the number of leapfrog steps can be tuned. In our PIGP implementation, we recommend fixing the number of leapfrog steps (set as $200$ in this paper), and tuning the leapfrog step size automatically during the burn-in period to achieve an acceptance rate between $60\%$ and $90\%$.
\subsection{Remarks on Optimization Algorithms for Maximizing Posterior Density}\label{spsubsec:map_opt_alg}
In this section, we describe the optimization algorithms for obtaining MAP from the posterior density (\ref{eq:map_opt_sig_e}).\par
As discussed in Section \ref{subsec:PIGPI1}, modularization is employed. 
We assign a simple two-step procedure to obtain the MAP $\hat{\bm{\theta}}, \hat{u}(\bm{I}), \hat{\sigma}_e^2$. 
In the first step, we train the GP model from the observation data, which can be found in Section \ref{spsubsec:train_GP}.
The hyper-parameters are obtained by optimizing the likelihood function of the GP.  
In the second step, we run an optimization algorithm to optimize (\ref{eq:map_opt_sig_e}) by choosing the initial guess for $\bm{\theta}$ (obtained from the two stage method) and initialized value for $u(\bm{I}),\sigma_e^2$ (obtained from the GP prediction and the MLE for $\sigma_e$) as starting point.
We recommend Adam provided in PyTorch in the second optimization stage for the full data case.
We recommend LBFGS provided in PyTorch in both optimization stages when there are censored components.\par
Note that we assume a compact support for $\sigma_e^2$, i.e., $\sigma_e^2\in [10^{-6}\psi_1, \psi_1]$.
In other words, we assume that the standard deviation of random error should not exceed the standard deviation of the GP model, i.e., the relative error is assumed less than 1. 
On the other hand, if the random error is negligible, we still add a nugget term ($10^{-6}\psi_1$) to avoid numerical issues.
\section{Supplementary Results for Examples in Section \ref{sec:simulation_res}}\label{spsec:simulation_complement}
In this section, we list the supplementary results for examples in Section  \ref{sec:simulation_res}.

\subsection{Supplementary Results for Example in Section \ref{subsec:egreaction_diffution}}\label{spsubsec:egreaction_diffution}
We list the supplementary results for example in Section \ref{subsec:egreaction_diffution}.
First, the bias and its RMSE for simulation results in Section \ref{subsubsec:eg_reaction_diffution_full} are listed in Table \ref{sptb:compare_map_multi_dim}.
The results shown in Table \ref{sptb:compare_map_multi_dim} are based on the simulation conducted in Section \ref{subsubsec:eg_reaction_diffution_full}.
\begin{table}[H]
\caption{The mean±SD of bias of MAPs obtained by PIGP, parameter estimation obtained by MLE, TSM, and BOM across 100 simulated datasets. The proposed method is emphasized with boldface font.}\label{sptb:compare_map_multi_dim}
\begin{center}
\begin{small}
\begin{tabular}{|l|l|rl|rl|rl|rl|}
\hline
                                                                                        & $n$   & \multicolumn{2}{l|}{30}          & \multicolumn{2}{l|}{60}         & \multicolumn{2}{l|}{120}       & \multicolumn{2}{l|}{240}       \\ \hline
\multirow{4}{*}{\begin{tabular}[c]{@{}l@{}}$\theta_1$\\ $\times 10^{-3}$\end{tabular}}  & BFM   & 0.14±           & 0.23           & 0.12±           & 0.16          & 0.09±          & 0.12          & 0.07±          & 0.10          \\ \cline{2-10} 
                                                                                        & \textbf{PIGP} & \textbf{29.01±} & \textbf{10.77} & \textbf{10.35±} & \textbf{3.10} & \textbf{8.30±} & \textbf{2.48} & \textbf{8.05±} & \textbf{1.80} \\ \cline{2-10} 
                                                                                        & BOM   & 955.36±         & 2241.81        & 1148.92±        & 2457.89       & 1052.56±       & 2502.33       & 679.03±        & 1871.82       \\ \cline{2-10} 
                                                                                        & TSM   & 294.21±         & 153.72         & 96.77±          & 45.79         & 28.48±         & 20.91         & 12.56±         & 6.08          \\ \hline
\multirow{4}{*}{\begin{tabular}[c]{@{}l@{}}$\theta_2$ \\ $\times 10^{-3}$\end{tabular}} & BFM   & 0.02±           & 0.55           & 0.01±           & 0.34          & 0.01±          & 0.23          & -0.01±         & 0.18          \\ \cline{2-10} 
                                                                                        & \textbf{PIGP} & \textbf{0.20±}  & \textbf{3.01}  & \textbf{0.22±}  & \textbf{1.58} & \textbf{0.28±} & \textbf{1.14} & \textbf{0.80±} & \textbf{0.79} \\ \cline{2-10} 
                                                                                        & BOM   & 15.78±          & 114.21         & 22.16±          & 104.38        & 19.51±         & 117.35        & 28.67±         & 109.20        \\ \cline{2-10} 
                                                                                        & TSM   & -811.22±        & 711.83         & -286.75±        & 867.00        & -137.39±       & 188.58        & -85.02±        & 137.97        \\ \hline
\multirow{4}{*}{\begin{tabular}[c]{@{}l@{}}$\theta_3$ \\ $\times 10^{-3}$\end{tabular}} & BFM   & 0.04±           & 0.50           & -0.03±          & 0.36          & 0.00±          & 0.26          & -0.01±         & 0.17          \\ \cline{2-10} 
                                                                                        & \textbf{PIGP} & \textbf{0.86±}  & \textbf{1.13}  & \textbf{0.89±}  & \textbf{0.94} & \textbf{0.96±} & \textbf{0.63} & \textbf{1.09±} & \textbf{0.39} \\ \cline{2-10} 
                                                                                        & BOM   & 24.08±          & 113.74         & 22.01±          & 128.83        & 37.02±         & 112.51        & 33.84±         & 107.09        \\ \cline{2-10} 
                                                                                        & TSM   & -397.63±        & 226.55         & -173.55±        & 218.44        & -124.18±       & 179.92        & -30.81±        & 50.97         \\ \hline
\end{tabular}
\end{small}
\end{center}
\end{table}

The bias and its RMSE for simulation results in Section \ref{subsubsec:eg_reaction_diffution_censor} are listed in Table \ref{sptb:compare_map_censored}.
The results shown in Table \ref{sptb:compare_map_censored} are based on the simulation conducted in Section \ref{subsubsec:eg_reaction_diffution_censor}.
\begin{table}[H]
\caption{The mean±SD of bias of MAPs obteind by PIGP, parameter estimation obtained by MLE and BOM across 100 simulated datasets. The proposed method is emphasized with bold face font.}\label{sptb:compare_map_censored}    
\begin{center}
\begin{small}
\begin{tabular}{|l|l|rl|rl|rl|rl|}
\hline
                                                                                       & $n$            & \multicolumn{2}{l|}{30}            & \multicolumn{2}{l|}{60}          & \multicolumn{2}{l|}{120}         & \multicolumn{2}{l|}{240}         \\ \hline
\multirow{3}{*}{\begin{tabular}[c]{@{}l@{}}$\theta_1$\\ $\times 10^{-3}$\end{tabular}} & BFM            & 0.12±            & 0.18            & 0.09±           & 0.14           & 0.08±           & 0.14           & 0.05±           & 0.08           \\ \cline{2-10} 
                                                                                       & \textbf{PIGP} & \textbf{-1.10±}  & \textbf{56.59}  & \textbf{-4.28±} & \textbf{4.37}  & \textbf{-4.23±} & \textbf{3.20}  & \textbf{-2.26±} & \textbf{2.96}  \\ \cline{2-10} 
                                                                                       & BOM            & 883.59±          & 2474.97         & 1007.14±        & 2455.64        & 940.61±         & 2276.38        & 1040.75±        & 2569.32        \\ \hline
\multirow{3}{*}{\begin{tabular}[c]{@{}l@{}}$\theta_2$\\ $\times 10^{-3}$\end{tabular}} & BFM            & -0.06±           & 0.75            & 0.04±           & 0.63           & 0.10±           & 0.43           & -0.01±          & 0.31           \\ \cline{2-10} 
                                                                                       & \textbf{PIGP} & \textbf{10.39±}  & \textbf{17.95}  & \textbf{1.89±}  & \textbf{3.16}  & \textbf{1.09±}  & \textbf{1.85}  & \textbf{1.79±}  & \textbf{1.37}  \\ \cline{2-10} 
                                                                                       & BOM            & 60.39±           & 311.81          & 45.88±          & 306.82         & -21.70±         & 301.04         & 92.04±          & 298.32         \\ \hline
\multirow{3}{*}{\begin{tabular}[c]{@{}l@{}}$\theta_3$\\ $\times 10^{-3}$\end{tabular}} & BFM            & -0.04±           & 1.29            & 0.08±           & 1.11           & 0.15±           & 0.78           & -0.05±          & 0.55           \\ \cline{2-10} 
                                                                                       & \textbf{PIGP} & \textbf{140.63±} & \textbf{180.26} & \textbf{29.04±} & \textbf{47.62} & \textbf{17.93±} & \textbf{28.22} & \textbf{26.04±} & \textbf{16.84} \\ \cline{2-10} 
                                                                                       & BOM            & 91.07±           & 558.94          & 64.97±          & 575.69         & -59.85±         & 575.64         & 154.57±         & 555.43         \\ \hline
\end{tabular}
\end{small}    
\end{center}
\end{table}

\end{document}